\input amstex
\documentstyle{amams} % input Annals of Mathematics macros.
\document
\def\bdar#1{\phantom{\scriptstyle #1}\Big\downarrow{\scriptstyle #1}}
\def\stackrel#1#2{\overset{#1}\to {#2}}
\def\navs{\noalign{\vskip5pt}}

\def\joinrel{\mathrel{\mkern-4mu}}
\def\relbar{\mathrel{\smash-}}
\def\lrar{\relbar\joinrel\relbar\joinrel\rightarrow}
\def\dlines{{\matrix
\relbar\joinrel\relbar\joinrel\relbar\joinrel\relbar\\ \noalign{\vskip-10pt}
\relbar\joinrel\relbar\joinrel\relbar\joinrel\relbar\endmatrix}}
\annalsline{158}{2003}
\received{May 22, 2001}
\startingpage{593}
\font\tenrm=cmr10

\catcode`\@=11
\font\twelvemsb=msbm10 scaled 1100

%\font\ninemsb=msbm7 scaled 1100%msbm9
\font\ninemsb=msbm10 scaled 800
\newfam\msbfam
\textfont\msbfam=\twelvemsb  \scriptfont\msbfam=\ninemsb
  \scriptscriptfont\msbfam=\ninemsb
\def\msb@{\hexnumber@\msbfam}
\def\Bbb{\relax\ifmmode\let\next\Bbb@\else
 \def\next{\errmessage{Use \string\Bbb\space only in math
mode}}\fi\next}
\def\Bbb@#1{{\Bbb@@{#1}}}
\def\Bbb@@#1{\fam\msbfam#1}
\catcode`\@=12

 \catcode`\@=11
\font\twelveeuf=eufm10 scaled 1100
\font\teneuf=eufm10
\font\nineeuf=eufm7 scaled 1100%eufm9
\newfam\euffam
\textfont\euffam=\twelveeuf  \scriptfont\euffam=\teneuf
  \scriptscriptfont\euffam=\nineeuf
\def\euf@{\hexnumber@\euffam}
\def\frak{\relax\ifmmode\let\next\frak@\else
 \def\next{\errmessage{Use \string\frak\space only in math
mode}}\fi\next}
\def\frak@#1{{\frak@@{#1}}}
\def\frak@@#1{\fam\euffam#1}
\catcode`\@=12

%--------------- Author macros ---------------
\def\scirc{\raise.2ex\hbox{${\scriptstyle\circ}$}}
\def\ssb{\raise.2ex\hbox{${\scriptscriptstyle\bullet}$}}
\def\asD{\langle D_{\ssb}\rangle}
\def\aD{\langle D\rangle}
\def\osX{\overline{X}_{\ssb}}
\def\cA{{\Cal A}}
\def\cC{{\Cal C}}
\def\cD{{\Cal D}}
\def\cH{{\Cal H}}
\def\cK{{\Cal K}}
\def\cM{{\Cal M}}
\def\cO{{\Cal O}}
\def\cP{{\Cal P}}
\def\cV{{\Cal V}}

\def\bC{{\Bbb C}}
\def\bG{{\Bbb G}}
\def\bP{{\Bbb P}}
\def\bQ{{\Bbb Q}}
\def\bZ{{\Bbb Z}}
\def\tD{\widetilde{D}}
\def\tG{\widetilde{G}}
\def\tM{\widetilde{M}}
\def\tX{\widetilde{X}}
\def\tY{\widetilde{Y}}
\def\tZ{\widetilde{Z}}
\def\tgG{\widetilde{\Gamma}}
\def\tcK{\widetilde{\Cal K}}
\def\hG{\widehat{G}}
\def\hM{\widehat{M}}
\def\hgG{\widehat{\Gamma}}
\def\oG{\overline{G}}
\def\ok{\overline{k}}
\def\oP{\overline{P}}
\def\oS{\overline{S}}
\def\oT{\overline{T}}
\def\oX{\overline{X}}
\def\oY{\overline{Y}}

\def\la{{\langle}}
\def\ra{{\rangle}}
\def\et{\text{\rm \'et}}
\def\DR{\text{\rm DR}}
\def\eff{\text{\rm eff}}
\def\lf{\text{\rm lf}}
\def\sa{\text{\rm sa}}
\def\nsqp{\text{\rm nsqp}}
\def\alg{\text{\rm alg}}
\def\Zar{\text{\rm Zar}}

\def\Lie{\hbox{{\rm Lie}}\,}
\def\Dec{\hbox{{\rm Dec}}}
\def\Pic{\hbox{{\rm Pic}}}
\def\div{\text{\rm div}}
\def\NS{\hbox{{\rm NS}}}
\def\CH{\hbox{{\rm CH}}}
\def\Spec{\hbox{{\rm Spec}}}
\def\Ext{\hbox{{\rm Ext}}}
\def\Gr{\text{{\rm Gr}}}
\def\Im{\hbox{{\rm Im}}}
\def\Ker{\hbox{{\rm Ker}}}
\def\Coker{\hbox{{\rm Coker}}}
\def\Hom{\hbox{{\rm Hom}}}
\def\MHS{\text{{\rm MHS}}}
\def\an{\text{\rm an}}
\def\AH{\text{\rm AH}}
\def\fr{\text{\rm fr}}
\def\tor{\text{\rm tor}}
\def\red{\text{\rm red}}
\def\simto{\buildrel\sim\over\to}

%-------------- Author entries --------------------
\font\eti=cmr10 scaled \magstep3
\title{Deligne's conjecture on {\eti 1}-motives} %Article title
\shorttitle{Deligne's conjecture on 1-motives} % Shortened version for headline title
% Acknowledgements: Please enter all acknowledgements here.

 \twoauthors{L.~Barbieri-Viale, A.~Rosenschon,}{M.~Saito}
   \institutions{University of Rome {\ninepoint\it La Sapienza\/}, Rome, Italy
\\
{\eightpoint {\it E-mail address\/}: Luca.Barbieri-Viale\@uniroma1.it}
\\
\vglue6pt
Duke University, Durham, NC\\
{\eightpoint {\it Curent address\/}}: University of Buffalo, The State University of New York, Buffalo, 
\phantom{{\eightpoint {\it Curent addressessiss:}}: }NY\\
{\eightpoint {\it E-mail address\/}: axr\@math.buffalo.edu}\\
\vglue6pt
RIMS Kyoto University, Kyoto, Japan\\
{\eightpoint {\it E-mail address\/}: msaito\@kurims.kyoto-u.ac.jp
}}

\vglue6pt
\centerline{\bf Abstract}
\vglue6pt
We reformulate a conjecture of Deligne on $1$-motives by using
the integral weight filtration of Gillet and Soul\'e on cohomology,
and prove it.
This implies the original conjecture up to isogeny.
If the degree of cohomology is at most two, we can prove the conjecture
for the Hodge realization without isogeny, and even for $1$-motives with
torsion.
 
\vfill
 
Let
$ X $ be a complex algebraic variety.
We denote by
$ {H}_{(1)}^{j}(X,\bZ) $ the maximal mixed Hodge
structure of type
$ \{(0,0) $,
$ (0,1) $,
$ (1,0) $,
$ (1,1)\} $ contained in
$ H^{j}(X,\bZ) $.
Let
$ {H}_{(1)}^{j}(X,\bZ)_{\fr} $ be the quotient of
$ {H}_{(1)}^{j}(X,\bZ) $ by the torsion subgroup.
P. Deligne ([10, 10.4.1]) conjectured that the $1$-motive corresponding to
$ {H}_{(1)}^{j}(X,\bZ)_{\fr} $ admits a purely algebraic description,
that is, there should exist a $1$-motive
$ M_{j}(X)_{\fr} $ which is defined without using the associated
analytic space, and whose image
$ r_{\cH}(M_{j}(X)_{\fr}) $ under the Hodge realization functor
$ r_{\cH} $ (see {\it loc.~cit}.\ and (1.5) below) is canonically
isomorphic to
$ {H}_{(1)}^{j}(X,\bZ)_{\fr}(1) $
(and similarly for the
$ l $-adic and de Rham realizations).

This conjecture has been proved for curves [10],
for the second cohomology of projective surfaces [9],
and for the first cohomology of any varieties [2]
(see also [25]).
In general, a careful analysis of the weight spectral sequence
in Hodge theory leads us to a candidate for
$ M_{j}(X)_{\fr} $ up to isogeny (see also [26]).
However, since the torsion part cannot be handled by Hodge theory,
it is a rather difficult problem to solve the conjecture without isogeny.

In this paper, we introduce the notion of an {\it effective} $1$-motive
which admits torsion.
By modifying morphisms, we can get an abelian
category of\break $1$-motives which admit torsion, and prove that this is
equivalent to the category of graded-polarizable mixed
$ \bZ $-Hodge structures of the above type.
However, our construction gives in general {\it nonreduced}
effective $1$-motives, that is, the discrete part has torsion and
its image in the semiabelian variety is nontrivial.
\vfill
\footnoterule
{\eightpoint 

\vglue-4pt
{\it Key words and phrases}. $1$-motive, weight filtration, Deligne cohomology,
Picard group.} \eject

Let
$ Y $ be a closed subvariety of
$ X $.
Using an appropriate `resolution', we can define a canonical integral
weight filtration
$ W $ on the relative cohomology
$ H^{j}(X,Y;\bZ) $.
This is due to Gillet and Soul\'e ([14, 3.1.2]) if
$ X $ is proper.
See also (2.3) below.
Let
$ {H}_{(1)}^{j}(X,Y;\bZ) $ be the maximal mixed Hodge structure
of the considered type contained in
$ H^{j}(X,Y;\bZ) $.
It has the induced weight filtration
$ W $, and so do its torsion part
$ {H}_{(1)}^{j}(X,Y;\bZ)_{\tor} $ and its free part
$ {H}_{(1)}^{j}(X,Y;\bZ)_{\fr} $.
Using the same resolution as above, we construct the desired effective
$1$-motive
$ M_{j}(X,Y) $.
In general, only its free part
$ M_{j}(X,Y)_{\fr} $ is independent of the choice of the resolution.
By a similar idea, we can construct the derived relative Picard groups
together with an exact sequence similar to Bloch's localization
sequence for higher Chow groups [7];  see (2.6).
Our first main result shows a close relation between the nonreduced
structure of our $1$-motive and the integral weight filtration:

\nonumproclaim{0.1.~Theorem} 
There exists a canonical isomorphism of mixed Hodge structures
$$
\phi_{\fr} : r_{\cH}(M_{j}(X,Y))_{\fr}(-1) \to
W_{2}{H}_{(1)}^{j}(X,Y;\bZ)_{\fr},
$$
such that the semiabelian part and the torus part of
$ M_{j}(X,Y) $ correspond respectively to
$ W_{1}{H}_{(1)}^{j}(X,Y;\bZ)_{\fr} $ and
$ W_{0}{H}_{(1)}^{j}(X,Y;\bZ)_{\fr} $.
A quotient of its discrete part by some torsion subgroup is isomorphic
to
$ {\Gr}_{2}^{W}{H}_{(1)}^{j}(X,Y;\bZ) $.
Furthermore{\rm ,} similar assertions hold for the
$ l $\/{\rm -}\/adic and de Rham realizations.
\endproclaim

This implies Deligne's conjecture for the relative cohomology up to
isogeny.
As a corollary, the conjecture without isogeny is reduced
to:
$$
{H}_{(1)}^{j}(X,Y;\bZ)_{\fr} = W_{2}{H}_{(1)}^{j}(X,Y;
\bZ)_{\fr}.
$$
This is satisfied if the
$ {\Gr}_{q}^{W}H^{j}(X,Y;\bZ) $ are torsion-free for
$ q > 2 $.
The problem here is that we cannot rule out the possibility of the
contribution of the torsion part of
$ {\Gr}_{q}^{W}H^{j}(X,Y;\bZ) $ to
$ {H}_{(1)}^{j}(X,Y;\bZ)_{\fr} $.
By construction,
$ M_{j}(X,Y) $ does not have information on
$ W_{1}{H}_{(1)}^{j}(X,Y;\bZ)_{\tor} $, and the morphism
$ \phi_{\fr} $ in (0.1) is actually induced by a morphism of mixed
Hodge structures
$$
\phi : r_{\cH}(M_{j}(X,Y))(-1) \to W_{2}{H}_{(1)}^{j}
(X,Y;\bZ)/W_{1}{H}_{(1)}^{j}(X,Y;\bZ)_{\tor}.
$$

\nonumproclaim{{0.2.~Theorem}}
The composition of
$ \phi $ and the natural inclusion
$$
r_{\cH}(M_{j}(X,Y))(-1) \to {H}_{(1)}^{j}
(X,Y;\bZ)/W_{1}{H}_{(1)}^{j}(X,Y;\bZ)_{\tor}
$$
is an isomorphism if
$ j \le 2 $ or if
$ j = 3 $,
$ X $ is proper{\rm ,} and has a resolution of singularities whose third
cohomology with integral coefficients is torsion\/{\rm -}\/free{\rm ,}
and whose second cohomology is of type
$ (1,1) $. \pagebreak
\endproclaim 

The proof of these theorems makes use of a cofiltration on a complex
of varieties, which approximates the weight filtration, and simplifies
many arguments.
The key point in the proof is the comparison of the extension classes
associated with a $1$-motive and a mixed Hodge structure, as indicated
in Carlson's paper [9].
This is also the point which is not very clear in [26].
We solve this problem by using the theory of mixed Hodge complexes
due to Deligne [10] and Beilinson [4].
For the comparison of algebraic structures on the Picard group,
we use the theory of admissible normal functions [29].
This also shows the representability of the Picard type functor.
However, for an algebraic construction of the semiabelian part of the
$1$-motive
$ M_{j}(X,Y) $, we have to verify the representability in a purely
algebraic way [2] (see also [26]).
The proof of (0.2) uses the weight spectral sequence [10] with integral
coefficients, which is associated to the above resolution; see (4.4).
It is then easy to show

\nonumproclaim{{0.3.~Proposition}}
Deligne\/{\rm '}\/s conjecture without isogeny is true if
$ {E}_{\infty}^{p,j-p} $ is torsion\/{\rm -}\/free for
$ p \le j - 3 $.
The morphism
$ \phi $ is injective if
$ {E}_{2}^{p,j-1-p} = 0 $ for\break\vglue-11pt\noindent 
$ p \le j - 4 $ and
$ {E}_{1}^{j-3,2} $ is of type
$ (1,1) $.
\endproclaim

The paper is organized as follows.
In Section~1 we review the theory of $1$-motives with torsion.
In Section~2, the existence of a canonical integral filtration is deduced
from [17] by using a complex of varieties.
(See also [14].)
In Section~3, we construct the desired $1$-motive by using a cofiltration
on a complex of varieties, and show the compatibility for the
$ l $-adic and de Rham realizations.
After reviewing mixed Hodge theory in Section~4, we prove the main
theorems in Section~5.

\demo{Acknowledgements} The first and second authors would like to thank
the European community Training and Mobility of Researchers Network titled\break
{\it Algebraic K\/{\rm -}\/Theory{\rm ,} Linear Algebraic Groups and Related Structures}
for\break financial support.
\enddemo

{\it Notation}.
In this paper, a variety means a separated reduced scheme of finite type
over a field.

\section{$1$-Motives}

We explain the theory of $1$-motives with torsion by modifying slightly
[10].
This would be known to some specialists.

\vglue12pt {1.1.} Let
$ k $ be a field of characteristic zero, and
$ \ok $ an algebraic closure of~$k$.
(The argument in the positive characteristic case is more complicated
due to the nonreduced part of finite commutative group schemes; see [22].)

An {\it effective} $1$-motive
$ M = [\Gamma \overset{f}\to\to G] $ over
$ k $ consists of a locally finite commutative group scheme
$ \Gamma /k $ and a semiabelian variety
$ G/k $ together with a morphism of
$ k $-group schemes
$ f : \Gamma \to G $ such that
$ \Gamma (\ok) $ is a finitely generated abelian group.
Note that
$ \Gamma $ is identified with
$ \Gamma (\ok) $ endowed with Galois action because
$ k $ is a perfect field.
Sometimes an effective $1$-motive is simply called a $1$-motive,
since the category of $1$-motives will be defined by modifying
only morphisms.
A locally finite commutative group scheme
$ \Gamma /k $ and a semiabelian variety
$ G/k $ are identified with $1$-motives
$ [\Gamma \to 0] $ and
$ [0 \to G] $ respectively.

An {\it effective} morphism of $1$-motives
$$
u = (u_{\lf}, u_{\sa}) : M = [\Gamma \overset{f}\to\to G] \to M'
= [\Gamma'\overset{f'}\to\to G']
$$
consists of morphisms of
$ k $-group schemes
$ u_{\lf} : \Gamma \to \Gamma' $ and
$ u_{\sa} : G \to G' $ forming a commutative diagram
(together with
$ f, f') $.
We will denote by
$$
\Hom_{\eff}(M,M')
$$
the abelian group of effective morphisms of $1$-motives.

An effective morphism
$ u = (u_{\lf}, u_{\sa}) $ is called {\it strict}, if the kernel
of
$ u_{\sa} $ is connected.
We say that
$ u $ is a {\it quasi-isomorphism} if
$ u_{\sa} $ is an isogeny and if we have a commutative diagram with
exact rows
$$
\matrix 
0 &\lrar& E &\lrar& \Gamma &\lrar& \Gamma' &\lrar& 0
\\
\noalign{\vskip6pt}
&&\Big\Vert&&\Big\downarrow&&\Big\downarrow\\
\\ \noalign{\vskip-6pt}
0 &\lrar& E &\lrar& G &\lrar& G' &\lrar& 0
\endmatrix
\leqno(1.1.1)
$$
(i.e.\  if the right half of the diagram is cartesian).

We define morphisms of $1$-motives by inverting quasi-isomorphisms
from the right; i.e.\  a morphism is represented by
$ u\scirc v^{-1} $ with
$ v $ a quasi-isomorphism.
More precisely, we define
$$
\Hom(M,M') = \varinjlim \Hom_{\eff}(\tM,M'),
\leqno(1.1.2)
$$
where the inductive limit is taken over isogenies
$ \tG \to G $, and
$ \tM = [\tgG \to \tG] $ with
$ \tgG = \Gamma \times_{G}\tG $.
(This is similar to the localization of a triangulated category in
[33].)
Here we may restrict to isogenies
$ n : G \to G $ for positive integers
$ n $, because they form a cofinal index subset.
Note that the transition
morphisms of the inductive system are injective by the surjectivity of
isogenies together with the property of fiber product.
By (1.2) below, $1$-motives form a category which will be denoted by
$ \cM_{1}(k) $.

Let
$ \Gamma_{\tor} $ denote the torsion part of
$ \Gamma $, and put
$ M_{\tor} = \Gamma_{\tor} \cap \Ker\, f $.
This is identified with
$ [M_{\tor} \to 0] $, and is called the torsion part of
$ M $.
We say that
$ M $ is {\it reduced} if \pagebreak
$ f(\Gamma_{\tor}) = 0 $, {\it torsion-free} if
$ M_{\tor} = 0 $, {\it free} if
$ \Gamma_{\tor} = 0 $, and {\it torsion} if
$ \Gamma $ is torsion and
$ G = 0 $ (i.e.\  if
$ M = M_{\tor}) $.
Note that
$ M $ is free if and only if it is reduced and torsion-free.
We say that
$ M $ has {\it split torsion}, if
$ M_{\tor} \subset \Gamma_{\tor} $ is a direct factor of
$ \Gamma_{\tor} $.

We define
$ M_{\fr} = [\Gamma /\Gamma_{\tor} \to G/f(\Gamma_{\tor})] $.
This is free, and is called the free part of
$ M $.
If
$ M $ is torsion-free,
$ M_{\fr} $ is naturally quasi-isomorphic to
$ M $.
This implies that
$ [\Gamma/M_{\tor} \to G] $ is quasi-isomorphic to
$ M_{\fr} $ in general, and (1.3) gives a short exact sequence
$$
0 \to M_{\tor} \to M \to M_{\fr}
\to 0.
$$

\vglue12pt {\it Remark}. If
$ M $ is free,
$ M $ is a $1$-motive in the sense of Deligne [10].
We can show
$$
\Hom_{\eff}(M,M') = \Hom(M,M')
\leqno(1.1.3)
$$
for
$ M, M' \in \cM_{1}(k) $ such that
$ M' $ is free.
This is verified by applying (1.1.1) to the isogenies
$ \tG \to G $ in (1.1.2).
In particular, the category of Deligne $1$-motives, denoted by
$ \cM_{1}(k)_{\fr} $, is a full subcategory of
$ \cM_{1}(k) $.
The functoriality of
$ M \mapsto M_{\fr} $ implies
$$
\Hom(M_{\fr},M') = \Hom(M,M')
\leqno(1.1.4)
$$
for
$ M \in \cM_{1}(k) $,
$ M' \in \cM_{1}(k)_{\fr} $.
In other words, the functor
$ M \mapsto M_{\fr} $ is left adjoint of the natural functor
$ \cM_{1}(k)_{\fr} \to \cM_{1}(k) $.

\nonumproclaim{{1.2.~Lemma}}
For any effective morphism
$ u : \tM \to M' $ and any quasi\/{\rm -}\/isomorphism
$ \tM' \to M' ${\rm ,} there exists a quasi\/{\rm -}\/isomorphism
$ \hM \to \tM $ together with a morphism
$ v : \hM \to \tM' $ forming a commutative
diagram.
Furthermore{\rm ,}
$ v $ is uniquely determined by the other morphisms and the
commutativity.
In particular{\rm ,} we have a well\/{\rm -}\/defined composition of morphisms of
$ 1 $\/{\rm -}\/motives {\rm (}\/as in {\rm [33])}
$$
\Hom(M,M') \times \Hom(M',M'') \to \Hom(M,M'').
\leqno(1.2.1)
$$
\endproclaim 

\demo{{P}roof}
For the existence of
$ \hM $, it is sufficient to consider the semiabelian part
$ \hG $ by the property of fiber product.
Then it is clear, because the isogeny
$ n : G' \to G' $ factors through
$ \tG' \to G' $ for some positive integer
$ n $, and it is enough to take
$ n : \tG \to \tG $.
We have the uniqueness of
$ v $ for
$ \hG $ since there is no nontrivial morphism of
$ \hG $ to the kernel of the isogeny
$ \tG' \to \tG $ which is a torsion group.
The assertion for
$ \hgG $ follows from the property of fiber product.
Then the first two assertions imply (1.2.1) using the injectivity of
the transition morphisms. \pagebreak
\enddemo

\nonumproclaim{{1.3.~Proposition}}
Let
$ u : M \to M' $ be an effective morphism of\break $1$\/{\rm -}\/motives.
Then there exists a quasi\/{\rm -}\/isomorphism
$ \tM' \to M' $ such that
$ u $ is lifted to a strict morphism
$ u' : M \to \tM' $ {\rm (}\/i.e.\ 
$ \Ker\, u'_{\sa} $ is connected\/{\rm ).}
In particular{\rm ,}
$ \cM_{1}(k) $ is an abelian category.
\endproclaim

{\it Proof}.
It is enough to show the following
assertion for the semiabelian variety part: There exists an isogeny
$ \tG' \to G' $ with a morphism
$ u'_{\sa} : G \to \tG' $ lifting
$ u_{\sa} $ such that
$ \Ker\,u'_{\sa} $ is connected.
(Indeed, the first assertion implies the existence of kernel and cokernel,
and their independence of the representative of a morphism is easy.)

For the proof of the assertion, we may assume that
$ \Ker\,u_{\sa} $ is torsion, dividing
$ G $ by the identity component of
$ \Ker\,u_{\sa} $.
Let
$ n $ be a positive integer annihilating
$ E := \Ker\,u_{\sa} $ (i.e.\ 
$ E \subset {}_{n}G $).
We have a commutative diagram
$$
\matrix
&&\hskip-14pt E &\stackrel{n}{\lrar}& E
\\
\navs
&&\hskip-14pt\Big\downarrow&&\Big\downarrow
\\
{}_{n}G &\stackrel{\iota}{\lrar}&\hskip-14pt G&\stackrel{n}{\lrar}& G
\\
\navs
\Big\downarrow&&\hskip-14pt\bdar{u_{\sa}}&&\bdar{u_{\sa}} 
\\
{}_{n}G'& \stackrel{\iota'}{\lrar}&\hskip-14pt G'&\stackrel{n}{\lrar}& G'.
\endmatrix
\leqno(1.3.1)
$$  
Let
$ \tG' $ be the quotient of
$ G' $ by
$ u_{\sa}\iota({}_{n}G) $, and let
$ q : G' \to \tG' $ denote the projection.
Since
$ u_{\sa}\iota({}_{n}G) \subset \iota'({}_{n}G') $, there is a canonical
morphism
$ q' :\break \tG' \to G' $ such that
$ q'q = n : G' \to G' $.
Then the
$ u_{\sa} $ in the right column of the diagram is lifted to a morphism
$ u'_{\sa} : G \to \tG' $ (whose composition with
$ q' $ coincides with
$ u_{\sa} $), because
$ G $ is identified with the quotient of
$ G $ by
$ {}_{n}G $.
Furthermore,
$ \Im\,u'_{\sa} $ is identified with the quotient of
$ G $ by
$ {}_{n}G + E $, and the last term coincides with
$ {}_{n}G $ by the assumption on
$ n $.
Thus
$ u'_{\sa} $ is injective, and the assertion follows.
\vglue12pt
 
{\it Remark}. An isogeny of semiabelian varieties
$ G' \to G $ with kernel
$ E $ corresponds to an {\it injective} morphism of $1$-motives
$$
[0 \to G'] \to [E \to G'] = [0 \to G].
$$
 
\nonumproclaim{{1.4.~Lemma}}
Assume
$ k $ is algebraically closed.
Then{\rm ,} for a $1$\/{\rm -}\/motive
$ M ${\rm ,} there exists a quasi\/{\rm -}\/isomorphism
$ M' \to M $ such that
$ M' = [\Gamma' \overset{f'}\to\to G'] $ has split torsion.  
\endproclaim 

\demo{{P}roof} Let
$ n $ be a positive integer such that
$ E := \Gamma_{\tor} \cap \Ker\,f $ is annihilated by
$ n $.
Then
$ G' $ is given by
$ G $ with isogeny
$ G' \to G $ defined by the multiplication by
$ n $.
Let
$ \Gamma' = \Gamma \times_{G} G' $.
We have a diagram of the nine lemma  
$$
\matrix
&&{}_{n}G &\dlines& {}_{n}G
\\
&&\Big\downarrow&&\Big\downarrow
\\ \navs
E &\lrar& \Gamma'_{\tor} &\lrar& f'(\Gamma'_{\tor})
\\
\Big\Vert&&\Big\downarrow&&\Big\downarrow
\\ \navs
E &\lrar& \Gamma_{\tor} &\lrar& f(\Gamma_{\tor}).
\endmatrix
\leqno(1.4.1)
$$
The
$ l $-primary torsion subgroup of
$ G $ is identified with the quotient of
$ V_{l}G := T_{l}G \otimes_{\bZ_{l}}\bQ_{l} $ by
$ M := T_{l}G $.
Let
$ M' $ be the
$ \bZ_{l} $-submodule of
$ V_{l}G $ such that
$ M' \supset M $ and
$ M'/M $ is isomorphic to the
$ l $-primary part of
$ f(\Gamma_{\tor}) $.
Then there exists a basis
$ \{e_{i}\}_{1\le i\le r} $ of
$ M' $ together with integers
$ c_{i} \,(1 \le i \le r) $ such that
$ \{l^{c_{i}}e_{i}\}_{1\le i\le r} $ is a basis of
$ M $.
So the assertion is reduced to the following,
because the assumption on the second exact sequence
$$
0 \to {}_{n}G \to f'(\Gamma'_{\tor}) \to
f(\Gamma_{\tor}) \to 0
$$
is verified by the above argument.

\nonumproclaim{Sublemma}
Let
$ 0 \to A_{i} \to B_{i} \to C \to 0 $ be
short exact sequences of finite abelian groups for
$ i = 1, 2 $.
Put
$ B = B_{1} \times_{C} B_{2} $.
Assume that the second exact sequence {\rm (}\/i.e.{\rm ,} for
$ i = 2) $ is the direct sum of
$$
0 \to \bZ/n\bZ \to \bZ /nb_{j}
\bZ \to \bZ/b_{j}\bZ \to 0,
$$
such that
$ A_{1} $ is annihilated by
$ n $.
Then the projection
$ B \to B_{2} $ splits.
\endproclaim
 
{\it Proof.}
We see that
$ B $ corresponds to
$ (e_{1},e_{2}) \in \Ext^{1}(C,A_{1} \times A_{2}) $, where the
$ e_{i} \in \Ext^{1}(C,A_{i}) $ are defined by the exact sequences.
Then it is enough to construct a morphism
$ u : A_{2} \to A_{1} $ such that
$ e_{1} $ is the composition of
$ e_{2} $ and
$ u $, because this implies an automorphism of
$ A_{1} \times A_{2} $ over
$ A_{2} $ which is defined by
$ (a_{1}, a_{2}) \mapsto (a_{1}-u(a_{2}), a_{2}) $ so that
$ (e_{1}, e_{2}) $ corresponds to
$ (0, e_{2}) $.
(Indeed, it induces an automorphism of
$ B $ over
$ B_{2} $ so that
$ e_{1} $ becomes
$ 0 $.)
But the existence of such
$ u $ is clear by hypothesis.
This completes the proof of (1.4).
\vglue12pt

The following is a generalization of Deligne's construction
([10, 10.1.3]).

\nonumproclaim{{1.5.~Proposition}}
If
$ k = \bC ${\rm ,} we have an equivalence of categories
$$
r_{\cH} : \cM_{1}(\bC) \simto \MHS_{1},
\leqno(1.5.1)
$$
where
$ \MHS_{1} $ is the category of mixed
$ \bZ $-Hodge structures
$ H $ of type
$$ \{(0,0)  ,
  (0,-1)  ,
  (-1,0)  ,
  (-1,-1)\} $$  such that
$ {\Gr}_{-1}^{W}H_{\bQ} $ is polarizable.
\endproclaim

\demo{{P}roof} The argument is essentially the same
as in [10].
For a $1$-motive
$ M = [\Gamma \overset{f}\to\to G] $, let
$ \Lie G \to G $ be the exponential map, and
$ \Gamma_{1} $ be its kernel.
Then we have a commutative diagram with exact rows
$$
\matrix
0 &\lrar&  \Gamma_{1} &\lrar&  H_{\bZ} &\lrar&  \Gamma &\lrar&  0
\\ \navs
&&\Big\Vert&&\Big\downarrow&&\Big\downarrow
\\ \navs
0 &\lrar&  \Gamma_{1} &\lrar&  \Lie G &\lrar&  G &\lrar&  0
\endmatrix
\leqno(1.5.2)
$$
which defines
$ H_{\bZ} $, and
$ F^{0}H_{\bC} $ is given by the kernel of the projection
$$
H_{\bC} := H_{\bZ} \otimes_{\bZ} \bC \to \Lie G.
$$
We get
$ W_{-1}H_{\bQ} $ from
$ \Gamma_{1} $, and
$ W_{-2}H_{\bQ} $ from the corresponding exact sequence for the torus
part of
$ G $.
(See also Remark below.)

We can verify that
$ H_{\bZ} $ and
$ F^{0} $ are independent of the representative of
$ M $ (i.e.\  a quasi-isomorphism induces isomorphisms of
$ H_{\bZ} $ and
$ F^{0}). $ Indeed, for an isogeny
$ M' \to M $, we have a commutative diagram with exact rows
$$
\matrix
0 &\lrar&  \Gamma'_{1} &\lrar&  \Lie G' &\lrar&  G' &\lrar&  0
\\
&&\Big\downarrow&&\Big\Vert&&\Big\downarrow
\\ \navs
0 &\lrar&  \Gamma_{1} &\lrar&  \Lie G &\lrar&  G &\lrar&  0
\endmatrix
\leqno(1.5.3)
$$
and the assertion follows by taking the base change by
$ \Gamma \to G $.
So we get the canonical functor (1.5.1).
We show that this is fully faithful and essentially surjective.
(To construct a quasi-inverse, we have to choose a splitting of the
torsion part of
$ H_{\bZ} $ for any
$ H \in \MHS_{1} $.)

For the proof of the essential surjectivity, we may assume that
$ H $ is either torsion-free or torsion.
Note that we may assume the same for $1$-motives by (1.4).
But for these
$ H $ we have a canonical quasi-inverse as in [10].
Indeed, if
$ H $ is torsion-free, we lift the weight filtration
$ W $ to
$ H_{\bZ} $ so that the
$ {\Gr}_{k}^{W}H_{\bZ} $ are torsion-free.
Then we put
$$
\Gamma = {\Gr}_{0}^{W}H_{\bZ},\quad G = J(W_{-1}H) \,(=
{\Ext}_{\MHS}^{1}(\bZ,W_{-1}H)),
$$
(see [8]), and
$ f : \Gamma \to G $ is given by the boundary map
$$
\Hom_{\MHS}(\bZ, {\Gr}_{0}^{W}H) \to {\Ext}_{\MHS}^{1}
(\bZ,W_{-1}H)
$$
associated with
$ 0 \to W_{-1}H \to H \to {\Gr}_{0}^{W}H
\to 0 $.
It is easy to see that this is a quasi-inverse.
The quasi-inverse for a torsion
$ H $ is the obvious one.

As a corollary, we have the full faithfulness of
$ r_{\cH} $ for free $1$-motives using (1.1.3).
So it remains to show that (1.5.1) induces
$$
\Hom(M, M') = \Hom(r_{\cH}(M), r_{\cH}(M'))
\leqno(1.5.4)
$$
when
$ M = [\Gamma \to G] $ is free and
$ M' $ is torsion.
Put
$ H = r_{\cH}(M) $.
We will identify both
$ M' $ and
$ r_{\cH}(M') $ with a finite abelian group
$ \Gamma' $.

Let
$ W_{-1}M = [0 \to G], {\Gr}_{0}^{W}M = [\Gamma \to 0] $.
Then we have a short exact sequence
$$
0 \to \Hom({\Gr}_{0}^{W}M, M') \to \Hom(M, M') \to
\Hom(W_{-1}M, M') \to 0,
$$
because
$ \Ext^{1}({\Gr}_{0}^{W}M, M') = \Ext^{1}(\Gamma, \Gamma') = 0 $.
Since we have the corresponding exact sequence for mixed Hodge structures
and the assertion for
$ {\Gr}_{0}^{W}M $ is clear, we may assume
$ M = W_{-1}M $, i.e.,
$ \Gamma = 0 $.

Let
$ T(G) $ denote the Tate module of
$ G $.
This is identified with the completion of
$ H_{\bZ} $ using (1.5.3).
Then
$$
\Hom(M,M') = \Hom(T(G),\Gamma') = \Hom(H_{\bZ},\Gamma'),
$$
and the assertion follows.
\enddemo

{\it Remark}. Let
$ T $ be the torus part of
$ G $.
Then we get in (1.5.2) the integral weight filtration
$ W $ on
$ H := r_{\cH}(M) $ by
$$
W_{-1}H_{\bZ} = \Gamma_{1},\quad
W_{-2}H_{\bZ} = \Gamma_{1} \cap \Lie T.
\leqno(1.5.5)
$$

\vglue-2pt
\section{Geometric resolution}
\vglue-4pt

Using the notion of a complex of varieties together with some
arguments from [17] (see also [14], [16]), we show the existence
of a canonical integral weight filtration on cohomology.

\vglue6pt {2.1.} Let
$ \cV_{k} $ denote the {\it additive} category of
$ k $-varieties, where a morphism
$ X' \to X'' $ is a (formal) finite
$ \bZ $-linear combination
$ \sum_{i} [f_{i}] $ with
$ f_{i} $ a morphism of connected component of
$ X' $ to
$ X'' $.
It is identified with a cycle on
$ X' \times_{k} X'' $ by taking the graph.
(This is similar to a construction in [14].)
We say that a morphism
$ \sum_{i} n_{i}[f_{i}] $ is proper, if each
$ f_{i} $ is.
The category of
$ k $-varieties in the usual sense is naturally viewed as a
subcategory of the above category.
For a
$ k $-variety
$ X $, we have similarly the additive category
$ \cV_{X} $ consisting of {\it proper}
$ k $-varieties over
$ X $, where the morphisms are assumed to be defined over
$ X $ in the above definition.

Since these are additive categories, we can define the categories of
complexes
$ \cC_{k}, \cC_{X}, $ and also the categories
$ \cK_{k}, \cK_{X} $ where morphisms are considered up to
homotopy as in [33].
We will denote an object of
$ \cC_{X}, \cK_{X} $ (or
$ \cC_{k}, \cK_{k}) $ by
$ (X_{\ssb},d) $, where
$ d : X_{j} \to X_{j-1} $ is the differential, and will be often
omitted to simplify the notation.
The structure morphism is denoted by
$ \pi : X_{\ssb} \to X $.
(This lower index of
$ X_{\ssb} $ is due to the fact that we consider only contravariant
functors from this category.)
For
$ i \in \bZ $, we define the shift of complex by
$ (X_{\ssb}[i])_{p} = X_{p+i} $.
We say that
$ Y_{\ssb} $ is a closed subcomplex of
$ X_{\ssb} $ if the
$ Y_{i} $ are closed subvarieties of
$ X_{i} $, and are stable by the morphisms appearing in the differential
of
$ X_{\ssb} $.

We will denote by
$ \cC_{X}^{b} $ the full subcategory of
$ \cC_{X} $ consisting of bounded complexes, and by
$ \cC_{{X}^{\nsqp}}^{b} $ the full subcategory of
$ \cC_{X}^{b} $ consisting of complexes of smooth
quasi-projective varieties.
(Here nsqp stands for nonsingular and quasiprojective.)
Let
$ D $ be a closed subvariety of
$ X $.
We denote by
$ \cC_{{X\aD}^{\nsqp}}^{b} $ the full subcategory
of
$ \cC_{X^{\nsqp}}^{b} $ consisting of
$ X_{\ssb} $ such that
$ D_{j} := \pi^{-1}(D) \cap X_{j} $ is locally either a connected
component or a divisor with simple normal crossings for any
$ j $.
Here simple means that the irreducible components of
$ D_{j} $ are smooth.
For an integer
$ j $, let
$ \cC_{X}^{b,\ge j} $ denote the full subcategory of
$ \cC_{X}^{b} $ consisting of complexes
$ X_{\ssb} $ such that
$ X_{i} = \emptyset $ for
$ i < j $, and similarly for
$ \cC_{{X}^{\nsqp}}^{b,\ge j} $,
$ \cC_{{X\aD}^{\nsqp}}^{b,\ge j} $.
Replacing
$ \cC $ with
$ \cK $, we define similarly
$ \cK_{{X}^{\nsqp}}^{b} $,
$ \cK_{{X\aD}^{\nsqp}}^{b} $, etc.

We say that
$ X_{\ssb} \in \cK_{X}^{b} $ is {\it strongly acyclic}
if there exist
$ X'_{\ssb} \in \cK_{X}^{b} $ isomorphic to
$ X_{\ssb} $ in
$ \cK_{X}^{b} $ and a finite filtration
$ G $ on
$ X'_{\ssb} $ such that the restriction of
$ G $ to each component
$ X'_{j} $ is given by direct factors, and for each integer
$ i $ there exists a birational proper morphism of
$ k $-varieties
$ g : Y' \to Y $ together with a closed subvariety
$ Z $ of
$ Y $ satisfying the following condition:
Letting
$ Z' = (Y'\times_{Y}Z)_{\red} $, the morphism
$ g : Y' \setminus Z' \to Y \setminus Z $ is an isomorphism and
the graded piece
$ {\Gr}_{i}^{G}X'_{\ssb} $ is isomorphic in
$ \cK_{X}^{b} $ to the single complex
associated to
$$
\matrix
Z' &\lrar&  Y'
\\ \navs
\Big\downarrow&&\Big\downarrow
\\ \navs
Z &\lrar&  Y
\endmatrix
\leqno(2.1.1)
$$
up to a shift of complex.
Clearly, this condition is stable by mapping cone.
We say that a morphism
$ X'_{\ssb} \to X_{\ssb} $ in
$ \cC_{X}^{b} $ or
$ \cK_{X}^{b} $ is a {\it strong
quasi-isomorphism} if its mapping cone is strongly acyclic in
$ \cK_{X}^{b} $.
This condition is stable by compositions, using the octahedral axiom of
the triangulated category.
Similarly, if
$ vu $ and
$ u $ or
$ v $ are strongly acyclic, then so is the remaining.
(It is not clear whether the strongly acyclic complexes form a thick
subcategory in the sense of Verdier.)

We say that a proper morphism of
$ k $-varieties
$ X' \to X $ has the {\it lifting property} if it induces a
surjective morphism
$$
X'(K) \to X(K)
$$
for any field
$ K $ (see [14]), or equivalently, if any irreducible subvariety of
$ X $ can be lifted birationally to
$ X' $.
We say that a morphism
$ u : X' \to X $ in
$ \cV_{k} $ has the {\it lifting property}, if for any connected
component
$ X_{i} $ of
$ X $,
there exists a connected component
$ X'_{i'} $ of
$ X' $ such that the restriction of
$ u $ to
$ X'_{i'} $ is given by a proper morphism
$$
f_{i} : X'_{i'} \to X_{i}
$$
with coefficient
$ \pm 1 $ and
$ f_{i} $ has the lifting property.
We say that a morphism
$ u : X' \to X $ in
$ \cV_{k} $ is {\it of birational type} if for any irreducible
component
$ X_{i} $ of
$ X $,
there exists uniquely a connected component
$ X'_{i'} $ of
$ X' $ such that the restriction of
$ u $ to
$ X'_{i'} $ is defined by a birational proper morphism
$$
f_{i} : X'_{i'} \to X_{i}
$$
with coefficient
$ \pm 1 $,
and this gives a bijection between the irreducible components of
$ X' $ and
$ X $.

For
$ X_{\ssb} \in \cC_{X}^{b} $, we say that a morphism
$ u : X'_{\ssb} \to X_{\ssb} $ of
$ \cC_{X}^{b} $ is a {\it quasi-projective resolution over}
$ X\aD $, if
$ X'_{\ssb} \in \cC_{{X\aD}^{\nsqp}}^{b} $ and
$ u $ is a strong quasi-isomorphism in
$ \cK_{X}^{b} $.
We say that
$ u $ is a quasi-projective resolution {\it of degree}
$ \ge j $ over
$ X\aD $, if furthermore
$ X'_{\ssb}, X_{\ssb} \in \cC_{X}^{b,\ge j} $ and
$ u : X'_{j} \to X_{j} $ is of birational type.
We denote by
$ \cK_{{X\aD}^{\nsqp}}^{b}(X_{\ssb}) $ the
category of quasi-projective resolutions
$ u : X'_{\ssb} \to X_{\ssb} $ over
$ X\aD $ (which are morphisms in
$ \cC_{X}^{b} $).
A morphism of
$ u $ to
$ v $ is a morphism
$ w $ of the source of
$ u $ to that of
$ v $ in
$ \cK_{X}^{b} $ such that
$ u = vw $ in
$ \cK_{X}^{b} $.
If
$ X_{\ssb} \in \cC_{X}^{b,\ge j} $, we define similarly
$ \cK_{{X\aD}^{\nsqp}}^{b,\ge j}(X_{\ssb}) $
by assuming further the condition on degree
$ \ge j $.
\vglue2pt 
For
$ X_{\ssb} \in \cC_{X}^{b,\ge j} $ and a closed subcomplex
$ Y_{\ssb} $, we say that
$$ u : (X'_{\ssb},Y'_{\ssb})\break \to (X_{\ssb},
Y_{\ssb})$$ is a {\it smooth quasi-projective modification of
degree}
$ \ge j $, if 
$ X'_{\ssb} \in \cC_{{X\aD}^{\nsqp}}^{b,\ge j} $,\break
$ Y'_{\ssb} =
(X'_{\ssb}\times_{Y_{\ssb}}X_{\ssb})_{\red}$,
$ u : X'_{\ssb} \to X_{\ssb} $ is a proper morphism
inducing an isomor-\break\vglue-11pt \noindent phism
$ X'_{\ssb} \setminus Y'_{\ssb} \to X_{\ssb}
\setminus Y_{\ssb} $,
and
$ u : X'_{j} \to X_{j} $ is of birational type.

\vglue8pt {\it Remarks}. (i)
A birational proper morphism
$ f : X' \to X $ has the lifting property.
Indeed, according to Hironaka [18], there exists a variety
$ X'' $ together with morphisms
$ g : X'' \to X $ and
$ h : X'' \to X' $,
such that
$ fh = g $ and
$ g $ is obtained by iterating blowing-ups with smooth centers.
(Here we may assume that the centers are smooth using Hironaka's theory
of resolution of singularities.)
This implies that a proper morphism has the lifting property if the
generic points of the irreducible components can be lifted.

\vglue4pt
(ii) For
$ X_{\ssb} \in \cC_{X}^{b,\ge j} $ and a closed subcomplex
$ Y_{\ssb} $ such that
$ \dim Y_{\ssb} < \dim X_{\ssb} $,
there exists a smooth quasi-projective modification
$ (X'_{\ssb},Y'_{\ssb}) \to
(X_{\ssb},Y_{\ssb}) $ of degree
$ \ge j $ by replacing
$ Y_{\ssb} $ with a larger subcomplex of the same dimension if
necessary.
This follows from [17, I, 2.6], except the birationality of
$ X'_{j} \to X_{j} $,
because there are connected components of
$ X'_{\ssb} $ which are not birational to irreducible components of
$ X_{\ssb} $.
Indeed, if we denote by
$ Z_{i,a} $ the images of the irreducible components of
$ X_{k} \,(k \le i) $ by morphisms to
$ X_{i} $ which are obtained by composing morphisms appearing in the
differential of
$ X_{\ssb} $, then the connected components of
$ X'_{i} $ are `sufficiently blown-up' resolutions of singularities of
$ Z_{i,a}$, and are defined by increasing induction on
$ i $, lifting the differential of
$ X_{\ssb} $ (see {\it loc.~cit}\/).
However, if
$ Z_{j,a} $ is a \pagebreak proper closed subvariety of some irreducible component
$ X_{j,b} $ of
$ X_{j} $, we may replace the resolution of
$ Z_{j,a} $ by its lifting to the resolution of
$ X_{j,b} $ using the lifting property, because the differential
$ X_{j} \to X_{j-1} $ is zero.

\nonumproclaim{{2.2.~Proposition}}
For any
$ X_{\ssb} \in \cC_{X}^{b,\ge j} ${\rm ,} there exists a
quasi-projective resolution
$ X'_{\ssb} \to X_{\ssb} $ of degree
$ \ge j $ over
$ X\aD ${\rm ,} and the category
$ \cK_{{X\aD}^{\nsqp}_{\phantom{|}}}^{b,\ge j}(X_{\ssb}) $\break is
weakly directed in the following sense\/{\rm :}
For any
$ u_{i} \in \cK_{{X\aD}^{\nsqp}}^{b,\ge j}
(X_{\ssb})$\break $(i = 1, 2) ${\rm ,} there exists
$ u_{3} \in \cK_{{X\aD}^{\nsqp}}^{b,\ge j}
(X_{\ssb}) $ together with morphisms
$ u_{3} \to u_{i} $ in
$ \cK_{{X\aD}^{\nsqp}}^{b,\ge j}(X_{\ssb}) $.
\endproclaim

\demo{{P}roof}
We first show that
$ \cK_{{X\aD}^{\nsqp}}^{b,\ge j}(X_{\ssb}) $ is
nonempty by induction on
$ n := \dim X_{\ssb} $.
There exists a smooth quasi-projective modification
$$(X'_{\ssb},Y'_{\ssb}) \to
(X_{\ssb},Y_{\ssb}) $$ of degree
$ \ge j $ as in the above Remark (ii).
Then we have a strong quasi-isomorphism
$$
C(Y'_{\ssb} \to X'_{\ssb} \oplus Y_{\ssb})
\to X_{\ssb}
$$
where the direct sum means the disjoint union.
So it is enough to show by induction that
$ \tX_{\ssb} := C(Y'_{\ssb} \to Y_{\ssb}) $
has a strong quasi-isomorphism
$$
\tZ_{\ssb} \to \tX_{\ssb}
\leqno(2.2.1)
$$
in
$ \cC_{X}^{b} $
such that
$ \tZ_{\ssb} \in
\cC_{{X\aD}^{\nsqp}}^{b,\ge j} $ and for any
irreducible component
$ \tZ_{j,i} $ of
$ \tZ_{j} $ the restriction of the differential to some
irreducible component
$ \tZ_{j+1,i} $ of
$ \tZ_{j+1} $ is given by an isomorphism onto
$ \tZ_{j,i} $ with coefficient
$ \pm 1 $, under the inductive hypothesis:
$$
\tX_{j+1} \to \tX_{j}\,\,\,
\text{has the lifting property.}
\leqno(2.2.2)
$$
Indeed, admitting this,
$ \tZ_{\ssb} $ is then isomorphic to the mapping cone of
$$
\tZ'_{\ssb} \to \mathbold{\oplus}_{i}
C(\pm id : \tZ_{j,i} \to \tZ_{j,i})[-j]
$$
with
$ \tZ'_{\ssb} \in
\cC_{{X\aD}^{\nsqp}}^{b,\ge j} $,
and the mapping cone of
$ \pm id $ is isomorphic to zero in
$ \cK_{X}^{b} $.

To show (2.2.1), we repeat the above argument with
$ X_{\ssb} $ replaced by
$ \tX_{\ssb} $,
and get a smooth quasi-projective modification
$ (\tX'_{\ssb},\tY'_{\ssb}) \to
(\tX_{\ssb},\tY_{\ssb}) $.
By the lifting property (2.2.2), we may assume that for any
irreducible component
$ \tX_{j,i} $ of
$ \tX_{j} $,
the corresponding irreducible component
$ \tX'_{j,i} $ of
$ \tX'_{j} $ has a morphism
$ f_{i} $ to
$ \tX_{j+1} $ such that the composition of
$ f_{i} $ and
$ d : \tX_{j+1} \to \tX_{j} $ is the
canonical morphism
$ \tX'_{j,i} \to \tX_{j,i} $ up to a sign.
If
$ \dim \tX'_{j,i} = \dim \tX_{\ssb} $,
then
$ f_{i} $ induces a birational morphism to
$ \Im\, f_{i} $ and we may assume that there exists an irreducible
component
$ \tX'_{j+1,i} $ such that the restriction of
$ d $ to
$ \tX'_{j+1,i} $ is given by the isomorphism
$ \tX'_{j+1,i} \to \tX'_{j,i} $ by replacing
$ \tX'_{\ssb} $ if necessary, because the differential
$ d $ of
$ \tX'_{\ssb} $ is defined by lifting
$ d $ of
$ \tX_{\ssb} $ (see {\it loc.~cit}\/).
Then we can modify the morphism
$ \tX'_{j+1} \to \tX_{j+1} $ by using
$ f_{i} $ for
$ \dim \tX_{j,i} < \dim \tX_{\ssb} $,
and replace
$ \tX'_{j} $ with the union of the maximal dimensional
components.
So we may assume that
$ \tX'_{j} $ is equidimensional, because the modified
$ \tX'_{\ssb} \to \tX_{\ssb} $ still
induces an isomorphism over the complement of
$ \tY_{\ssb} $ by replacing
$ \tY_{\ssb} $ if necessary.
Here we may assume also that
$ \tY_{j+1} \to \tY_{j} $ has the lifting
property by taking
$ \tY_{\ssb} $ appropriately (due to (2.2.2) and the above
Remark (i)).
Then, considering the mapping cone of
$ \tY'_{\ssb} \to \tY_{\ssb} $,
the first assertion follows by induction.

The proof of the second assertion is similar.
Consider the shifted mapping cone (i.e.\  the first term has degree zero):
$$
X'_{\ssb} = [X_{1,\ssb}\oplus X_{2,\ssb} \to
X_{\ssb}],
\leqno(2.2.3)
$$
where the morphism is given by
$ u_{1} - u_{2} $.
Then
$ X'_{\ssb} \to X_{a,\ssb} $ is a strong
quasi-isomorphism.
Note that the composition of the canonical morphism
$ X'_{\ssb} \to X_{a,\ssb} $ and
$ u_{a} $ is independent of
$ a $ up to homotopy.

By definition, for any irreducible component
$ X'_{j-1,i} = X_{j,i} $ of
$ X'_{j-1} = X_{j} $,
there exist two connected components
$ Z_{i}, Z'_{i} $ of
$ X'_{j} $ such that the restrictions of
$ d $ to
$ Z_{i}, Z'_{i} $ are given by proper morphisms
$ Z_{i} \to X'_{j-1,i}, Z'_{i} \to X'_{j-1,i} $ which
have the lifting property (with coefficient
$ \pm 1) $.
Then by the same argument as above, we have a smooth quasi-projective
modification
$ u' : (X''_{\ssb},Y''_{\ssb}) \to (X'_{\ssb},
Y'_{\ssb}) $ of degree
$ \ge j - 1 $.
Here we may assume that the connected component
$ X''_{j-1,i} $ of
$ X''_{j-1} $ which is birational to
$ X'_{j-1,i} $ has morphisms to
$ Z_{i}, Z'_{i} $ factorizing the morphisms to
$ X'_{j-1,i} $, and
$ X''_{j} $ has two connected components such that the restriction of
$ d $ to each of these components is given by an isomorphism onto
$ X''_{j-1,i} $ (with coefficients
$ \pm 1) $.
We may also assume that
$ Y'_{j} \to Y'_{j-1} $ has the lifting property as before.

Then, applying the same argument to
$ C(Y''_{\ssb} \to Y'_{\ssb}) $,
and using induction on dimension, we get a strong quasi-isomorphism
$$
\tX_{\ssb} \to X'_{\ssb}
$$
such that
$ \tX_{\ssb} \in
\cC_{{X\aD}^{\nsqp}}^{b,\ge j-1} $,
and for any irreducible component
$ \tX_{j-1,i} $ of
$ \tX_{j-1} $,
$ \tX_{j} $ has two connected components such that the
restrictions of
$ d $ (resp.\ of the morphism to
$ X'_{j}) $ to these components are given by isomorphisms onto
$ \tX_{j-1,i} $ (resp.\ by birational proper morphisms to
$ Z_{i} $,
$ Z'_{i} $) with coefficients
$ \pm 1 $.
Thus
$ \tX_{\ssb}^{\phantom{|}}$ is isomorphic to the mapping cone of  
$$
\tX'_{\ssb} \to \mathbold{\oplus}_{i} C(\pm id :
\tX_{j-1,i} \to \tX_{j-1,i})[-j+1]
\leqno(2.2.4)
$$
where
$ \tX'_{\ssb} \in
\cC_{{X\aD}^{\nsqp}}^{b,\ge j-1} $.
So the second assertion follows.
\pagebreak

{\it  Remark.} It is not clear if for any
$ u_{i} \in \cK_{{X\aD}^{\nsqp}}^{b,\ge j}
(X_{\ssb})\,(i = 1, 2) $ and
$ v_{a} : u_{2} \to u_{1} $, there exists
$ u_{3} \in \cK_{{X\aD}^{\nsqp}}^{b,\ge j}
(X_{\ssb}) $ together with
$ w : u_{3} \to u_{2} $ such that
$ v_{1}^{\phantom{|}}w = v_{2}w $.
This condition is necessary to define an inductive limit over the
category
$ \cK_{{X\aD}^{\nsqp}}^{b,\ge j}(X_{\ssb}^{\phantom{|}}) $.
If we drop the condition on the degree
$ \ge j $, it can be proved for
$ \cK_{{X\aD}^{\nsqp}}^{b}(X_{\ssb}^{\phantom{|}}) $ by using the mapping cone
(2.2.3).
Indeed, let
$ K_{i,\ssb} $ be the source of
$ u_{i} $ for
$ i = 1, 2 $, and
$ K_{3,\ssb} $ the mapping cone of\break
$ (v_{1} - v_{2},0) : C(K_{2,\ssb} \to 0) \to C(K_{1,\ssb} \to
K_{\ssb})[-2] $, choosing a homotopy
$ h $ such that
$ dh + hd = u_{1}v_{1} - u_{1}v_{2} $.
Then
$ w : K_{3,\ssb} \to K_{2,\ssb}$ is given by the projection, and
$ v_{1}w - v_{2}w : K_{3,\ssb} \to K_{1,\ssb} $ factors through a
morphism
$ (v_{1} - v_{2},0) : K_{3,\ssb} \to C(K_{1,\ssb} \to K_{\ssb})[-1], $
which is homotopic to zero.

\nonumproclaim{{2.3.~Corollary}}
For a complex algebraic variety
$ X $ and a closed subvariety
$ Y ${\rm ,} there is a canonical integral weight filtration
$ W $ on the relative cohomology
$ H^{j}(X,Y;\bZ) $.
Furthermore{\rm ,} it is defined by a quasi\/{\rm -}\/projective resolution
$ \osX \to C(\oY \to \oX) $
of degree
$ \ge 0 $ over
$ \oX\aD ${\rm ,} where
$ \oX $ is a compactification of
$ X ${\rm ,}
$ \oY $ is the closure of
$ Y $ in
$ \oX $, and
$ D = \oX \setminus X $.
\endproclaim

{\it  Remarks.} (i) The first assertion is due to Gillet and Soul\'e
([14, 3.1.2]) in the case
$ X $ is proper (replacing
$ X_{\ssb} $ with a simplicial resolution).
It is expected that their integral weight filtration coincides with
ours.

(ii) If
$ X $ is proper, we have
$$
W_{i-1}H^{i}(X,Y;\bZ) = \Ker(H^{i}(X,Y;\bZ) \to
H^{i}(X',\bZ))
\leqno(2.3.1)
$$
for any resolution of singularities
$ X' \to X $ (see also {\it loc.\ cit.}).
Note that
$ \pi^{*} : H^{i}(X',\bZ) \to H^{i}(X'',\bZ) $
is injective for any birational proper morphism of smooth varieties
$ \pi : X'' \to X' $.

\vglue4pt {\it Proof of} (2.3).
The canonical mixed Hodge structure on the relative cohomology
can also be defined by using any quasi-projective resolution
$ \osX \to C(\oY \to \oX) $ as in [10].
This gives an integral weight filtration together with an integral
weight spectral sequence
$$
{E}_{1}^{p,q} = \mathbold{\oplus}_{k\ge 0} H^{q-2k}({\tD}_{p+k}^{k},\bZ)(-k)
\Rightarrow H^{p+q}(X_{\ssb},\bZ) = H^{p+q}(X,Y;\bZ),
\leqno(2.3.2)
$$
where
$ {\tD}_{j}^{k} $ the disjoint union of the intersections of
$ k $ irreducible components of~$ D_{j} $, and the cohomology is defined by taking the canonical
flasque resolution of Godement in the analytic or Zariski topology.
By (2.2) we get a set of integral weight filtrations on
$ H^{j}(X,Y;\bZ) $ which is directed with respect to the natural
ordering by the inclusion relation.
Then this is stationary \pagebreak by the noetherian property.
(It is constant if
$ X $ is proper; see (2.5) below.)
By the proof of (2.2) the limit is independent of the choice of
the compactification $ \oX $.
So the assertion follows.

\vglue12pt 2.4.\ {\it Definition.} For a complex of
$ k $-varieties
$ X_{\ssb} $ (see (2.1)), we define
$$
\Pic(X_{\ssb}) = H^{1}(X_{\ssb},
\cO_{{X}_{\ssb}}^{*})\,\,\, \text{(see [2])} .
\leqno(2.4.1)
$$
The right-hand side is defined by taking the canonical flasque
resolution of Godement which is compatible with the pull-back
by the differential of
$ X_{\ssb} $.
For a
$ k $-variety
$ X $ and closed subvariety
$ Y $, we define the {\it derived relative Picard groups} by
$$
{\bf Pic}(X,Y;i) = \varinjlim \Pic(X_{\ssb}[i]),
\leqno(2.4.2)
$$
where the inductive limit is taken over
$ X_{\ssb} \in \cK_{{X}^{\nsqp}}^{b}(C(Y \to X)) $.
If
$ Y $ is empty,
$ {\bf Pic}(X,Y;i) $ will be denoted by
$ {\bf Pic}(X,i) $, and
$ i $ will be omitted if
$ i = 0 $.

\vglue12pt {\it Remark.}
We can define similarly the derived relative Chow cohomology group by
$$
{\bf CH}^{p}(X,Y;i) = \varinjlim H^{p+i}(X_{\ssb},\cK_{p}),
\leqno(2.4.3)
$$
where
$ \cK_{p} $ is the Zariski sheafification of Quillen's higher
$ K $-group.
(In the case
$ X $ is smooth proper and
$ Y $ is empty, this is related to Bloch's higher Chow group for
$ i = 0, -1 $.)
 \vglue12pt

The following is a variant of a result of Gillet and Soul\'e
[14, 3.1], and gives a positive answer to the question in [2, 4.4.4].

\nonumproclaim{{2.5.~Proposition}}
Assume
$ X $,
$ Y $ proper.
Then a strong quasi\/{\rm -}\/isomor\-phism
$ u : X''_{\ssb} \to X'_{\ssb} $ in
$ \cC_{X^{\nsqp}}^{b} $ induces {\rm (}filtered{\rm )} isomorphisms
$$
u^{*} : \Pic(X'_{\ssb}[i]) \to \Pic(X''_{\ssb}[i]),\quad
u^{*} : (H^{i}(X'_{\ssb},\bZ),W) \to (H^{i}(X''_{\ssb},\bZ),W),
$$
where we assume
$ k = \bC $ for the second morphism.
In particular{\rm ,} the inductive system in {\rm (2.4.2)} is a constant
system in this case.
\endproclaim

\demo{{P}roof}
It is sufficient to show that
$ \Pic(X_{\ssb}[i]) = 0 $ and the
$ E_{1} $-complex
$ {}_{\bZ}E_{1}^{\ssb,q} $ of the integral weight spectral sequence
is acyclic,
if
$ X_{\ssb} $ is strongly acyclic and the
$ X_{j} $ are smooth.
(Note that the
$ E_{1} $-complex for
$ W $ is compatible with the mapping cone.)
Considering the
$ E_{1} $-complex of the spectral sequence
$$
{}_{P}{E}_{1}^{p,q} = H^{q}(X_{p},\cO_{{X}_{p}}^{*}) \Rightarrow
H^{p+q}(X_{\ssb},\cO_{{X}_{\ssb}}^{*}),
\leqno(2.5.1)
$$
it is enough to show the acyclicity of the complexes
$ {}_{P}E_{1}^{\ssb,q} $ and
$ {}_{\bZ}E_{1}^{\ssb,q} $ (where
$ {}_{P}E_{1}^{\ssb,q} = 0 $ for
$ q > 1 $; see (2.5.3) below).

By Gillet and Soul\'e ([14, 1.2]) this is further reduced to the
acyclicity of the Gersten complex of
$ X_{\ssb} \times V $ for any smooth proper variety
$ V $ because it implies that the image of the complex
$ X_{\ssb} $ in the category of complexes of varieties whose
differentials and morphisms are given by correspondences is
homotopic to zero.
Since the functor associating the Gersten complex preserves homotopy,
it is sufficient to show that the Gersten complex of
$$
Z' \to Z \oplus Y' \to Y
$$
is acyclic in the notation of (2.1.1)
(replacing it by the product with
$ V $).
Consider the subcomplex of the Gersten complex given by the points of
$ Z' $,
$ Z $,
$ Y' $,
$ Y $ contained in
$ Z' $,
$ Z $,
$ Z' $,
$ Z $ respectively.
It is clearly acyclic, and so is its quotient complex.
This shows the desired assertion.
(A similar argument works also for (2.4.3).)
\enddemo

{\it Remark.} Let
$ X $ be a smooth irreducible
$ k $-variety, and
$ k(X) $ the function field of
$ X $.
For closed subvariety
$ D $, let
$ k(X)_{X}^{*} $ and
$ \bZ_{D} $ denote the constant sheaf in the Zariski topology on
$ X $ and
$ D $ with stalk
$ k(X)^{*} $ and
$ \bZ $ respectively.
Then we have a flasque resolution
$$
0 \to \cO_{X}^{*} \to k(X)_{X}^{*}
\overset{\div}\to\to \mathbold{\oplus}_{D} \bZ_{D} \to 0,
\leqno(2.5.2)
$$
where the direct sum is taken over irreducible divisors
$ D $ on
$ X $.
In particular, we get
$$
\align
H^{i}(X,\cO_{X}^{*})
&=  0\quad \text{for }i > 1,
\tag {2.5.3}
\\
R^{i}j_{*}\cO_{X}^{*}
&=  0\quad \text{for }i > 0,
\tag {2.5.4}
\endalign
$$
for an open immersion
$ j : U \to X $.

\nonumproclaim{{2.6.~Proposition}}
There is a canonical long exact sequence
$$
\to {\bf Pic}(X,Y;i) \to
{\bf Pic}(X,i) \to
{\bf Pic}(Y,i) \to
{\bf Pic}(X,Y;i+1) \to .
\leqno(2.6.1)
$$
\endproclaim

{\it Remark.} This is an analogue of the localization sequence for
higher Chow groups [7].

\demo{{P}roof of {\rm (2.6)}}
The long exact sequence is induced by the distinguished triangle
$$
\to Y \overset{i}\to\to X \to C(Y \to X)
\to
$$
in
$ \cK_{X}^{b} $, because for any
quasi-projective resolutions
$ u : X_{\ssb} \to X $ and\break
$ v : Y_{\ssb} \to Y $, there exists
quasi-isomorphisms
$ u' : X'_{\ssb} \to X_{\ssb} $ and
$ v' : Y'_{\ssb} \to Y_{\ssb} $ together with
$ i' : Y'_{\ssb} \to X'_{\ssb} $ such that
$ u\scirc u'\scirc i' = i\scirc v\scirc v' $ in
$ \cK_{X}^{b} $ by using the mapping cone (2.2.3).
\enddemo

 2.7. {\it Remark}.
Assume
$ k = \bC $ and
$ X $ is proper.
Let
$ {H}_{\cD}^{i+2}(X,Y;\bZ(1)) $ denote the relative Deligne
cohomology.
See [3] and also (5.2) below.
Then we can show
$$
{\bf Pic}(X,Y;i) = {H}_{\cD}^{i+2}(X,Y;\bZ(1))
\quad \text{for }i \le 0.
\leqno(2.7.1)
$$
This is analogous to the canonical isomorphisms for
$$
\CH^{1}(X,i) = {H}_{2n-2+i}^{\AH}(X,\bZ(n-1))
$$
which holds for
$ i > 0 $ and any variety
$ X $ of dimension
$ n $ [31].

Assume
$ X $ is proper and normal.
Let
$ \cH_{X}^{1}\bZ(1) $ be the Zariski sheaf associated with
the presheaf
$ U \mapsto H^{1}(U,\bZ(1)) $.
By the Leray spectral sequence we get a natural injective morphism
$ {H}_{\Zar}^{1}(X,\cH_{X}^{1}\bZ(1)) \to
H^{2}(X,\bZ(1)) $.
See [1], [6].
Define a subgroup
$ {H}_{\cD}^{2}(X,\bZ(1))_{\alg} $ of the Deligne cohomology
$ {H}_{\cD}^{2}(X,\bZ(1)) $ (see (4.2) below) by the cartesian
diagram
$$
\matrix
{H}_{\cD}^{2}(X,\bZ(1))_{\alg} &\lrar& 
F^{1} \cap {H}_{\Zar}^{1}(X,\cH_{X}^{1}\bZ(1))
\\ \navs
\Big\downarrow&&\Big\downarrow
\\ \navs
{H}_{\cD}^{2}(X,\bZ(1)) &\lrar&  F^{1} \cap H^{2}(X,\bZ(1)).
\endmatrix
\leqno(2.7.2)
$$
Let
$ X_{\ssb} \to X $ be a quasi-projective resolution.
Then
$$
{H}_{\cD}^{2}(X,\bZ(1))_{\alg} =
\Im(\Pic(X) \to \Pic(X_{\ssb}) =
{H}_{\cD}^{2}(X,\bZ(1))).
\leqno(2.7.3)
$$
Indeed, if we put
$ \NS(X_{\ssb}) := \Im(\Pic(X_{\ssb}) \to H^{2}(X,\bZ(1))) $, and
similarly for
$ \NS(X) $, this follows from a result of Biswas and
Srinivas [6]:
$$
\NS(X) = F^{1} \cap {H}_{\Zar}^{1}(X,\cH_{X}^{1}\bZ(1)),
\leqno(2.7.4)
$$
by using
$$
\Coker(\Pic(X) \to \Pic(X_{\ssb})) =
\NS(X_{\ssb})/\NS(X).
\leqno(2.7.5)
$$
The last isomorphism follows from the morphism of long exact
sequences
$$
\matrix
H^{1}(X,\bZ(1)) &\lrar&  H^{1}(X,\cO_{X}) &\lrar& 
\Pic(X) &\lrar&  H^{2}(X,\bZ(1))
\\ \navs
\Big\Vert&&\bdar{(\ast)}&&\Big\downarrow&&\Big\Vert
\\ \navs
H^{1}(X_{\ssb},\bZ(1)) &\lrar& 
H^{1}(X_{\ssb},\cO_{X_{\ssb}}) &\lrar& 
\Pic(X_{\ssb}) &\lrar&  H^{2}(X_{\ssb},\bZ(1))
\endmatrix
$$
because
$ (*) $ is surjective by Hodge theory [10] (considering
the canonical morphisms of
$ H^{1}(X,\bC) $ to the source and the target of \pagebreak
$ (*) $). 

\section{Construction}

We construct $1$-motives associated with a complex of varieties,
and show the compatibility for the
$ l $-adic and de Rham realizations.
We assume
$ k $ is an algebraically closed field of characteristic zero.

\vglue12pt {3.1.} With the notation of (2.1), let
$ X_{\ssb} \in \cC_{k} $ be a complex of smooth
$ k $-varieties (see (2.1)), and
$ \osX $ a smooth compactification of
$ X_{\ssb} $ such that
$ D_{p} := \oX_{p} \setminus X_{p} $ is a divisor with simple
normal crossings.
We assume
$ X_{\ssb} $ is bounded below.
The reader can also assume that
$ (\osX,D_{\ssb}) $ is the mapping cone
$ (\oY_{\ssb},D'_{\ssb}) \to (\oX
_{\ssb},D_{\ssb}) $ of simplicial resolutions of
$ f : (\oY,D') \to (\oX,D) $ (see [10]), where
$ \oX, \oY $ are proper
$ k $-varieties with closed subvarieties
$ D, D' $ such that
$ D' = f^{-1}(D) $.

Let
$ j : X_{\ssb} \to \osX $ denote the
inclusion.
Put
$$
\cK = \Gamma (\osX,\cC^{\ssb}
(j_{*}\cO_{{X}_{\ssb}}^{*})),
$$
where
$ \cC^{\ssb} $ is the canonical flasque resolution of Godement
in the Zariski\break topology.

We define a {\it cofiltration}
$$
{}^{o}W^{\prime p}X_{\ssb}
$$
to be the quotient complexes of
$ X_{\ssb} $ consisting of
$ X_{i} $ for
$ i \ge p $ (and empty otherwise).
This is similar to the filtration ``b\^ete''
$ \sigma $ in [10].
It induces a decreasing filtration
$ W' $ on
$ \cK $ such that
$$
W^{\prime i}\cK = \Gamma ({}^{o}W^{\prime i}
\osX,\cC^{\ssb}
(j_{*}\cO_{{X}_{\ssb}}^{*}))).
$$

We define a cofiltration
$ {}^{o}W^{\prime \prime j}X_{\ssb} $ for
$ j = -1, 0, 1 $ by
$$
{}^{o}W^{\prime \prime -1}X_{\ssb} = X_{\ssb},\quad
{}^{o}W^{\prime \prime 0}X_{\ssb}
=\osX,\quad {}^{o}W^{\prime \prime 1}X_{\ssb}
= \emptyset .
$$
Since this depends on the compactification
$ \osX $, it is also denoted by
$$ {}^{o}W^{\prime \prime j}
\osX\asD $$
so that
$$
{}^{o}W^{\prime \prime -1}
\osX\asD
= \osX\asD,\quad
{}^{o}W^{\prime \prime 0}
\osX\asD
= \osX,\quad
{}^{o}W^{\prime \prime 1}
\osX\asD
= \emptyset .
$$
This corresponds to a decreasing filtration
$ W'' $ on
$ \cK $ such that
$$
W^{\prime \prime -1}\cK = \cK,\quad W^{\prime \prime 0}\cK
= \Gamma (\osX,
\cC^{\ssb}(\cO_{{\oX}_{\ssb}}^{*})),\quad
W^{\prime \prime 1}\cK =0.
$$
Then
$ {\Gr}_{W''}^{-1}\cK $ is quasi-isomorphic to
$ \Gamma (\tD_{\ssb},\bZ) $, and
$ \Gamma (\tD_{p},\bZ) \simeq \bZ^{\oplus r_{p}} $,
where
$ \tD_{\ssb} $ is the normalization of
$ D_{\ssb} $, and
$ r_{p} $ is the number of irreducible components of
$ D_{p} $.
(Indeed, a constant sheaf on an irreducible variety is flasque in
the Zariski topology.)

Let
$ W $ be the convolution of
$ W' $ and
$ W'' $ (i.e.,
$ W^{r} = \sum_{i+j=r} W^{\prime i} \cap W^{\prime \prime j}) $ so that
$$
{\Gr}_{W}^{r}\cK = \mathbold{\oplus}_{i+j=r} {\Gr}_{W'}^{i}{\Gr}_{W''}^{j}
\cK,
$$
(see [5, 3.1.2]).
This corresponds to the cofiltration
$ {}^{o}W^{r} $ such that
$$
{}^{o}W^{r}X_{\ssb}\,\,\, \text{(or }{}^{o}W^{r}
\osX\asD)
$$
consists of
$ X_{i} $ (or
$ \osX\asD $) for
$ i > r $, and
$ \oX_{r} $ for
$ i = r $.
Then we have a natural quasi-isomorphism
$$
\Gamma (\oX_{r}, \cC^{\ssb}
(\cO_{{\oX}_{r}}^{*}))[-r] \oplus
\Gamma (\tD_{r+1},\bZ)[-r-1] \to
{\Gr}_{W}^{r}\cK. 
\leqno(3.1.1) 
$$
Hence
$ H^{j}{\Gr}_{W}^{r}\cK $ vanishes unless
$ j = r $ or
$ r + 1 $, and
$$
H^{r}{\Gr}_{W}^{r}\cK = \Gamma (\oX_{r},
\cO_{{\oX}_{r}}^{*}),\quad H^{r+1}{\Gr}_{W}^{r} \cK =
\Pic(\oX_{r}) \oplus \Gamma (\tD_{r+1},\bZ),
\leqno(3.1.2)
$$
where
$ \Pic(\oX_{r}) = H^{1}(\oX_{r},
\cO_{{\oX}_{r}}^{*}) $ is the Picard group of
$ \oX_{r} $.

Let
$ \cK' = {\Gr}_{W''}^{0}\cK \,(= \Gamma (\osX,
\cC^{\ssb}(\cO_{{\oX}_{\ssb}}^{*}))) $ with
the induced filtration
$ W $.
Then we have the spectral sequence
$$
{E}_{1}^{p,q}=
\cases
H^{q}(\oX_{p},\cO_{{\oX}_{p}}^{*})
&\text{for $ p \ge r $}
\\
0 &\text{for $ p < r $}
\endcases
\Rightarrow H^{p+q}(W^{r}\cK'),
\leqno(3.1.3)
$$
which degenerates at
$ E_{3} $.
We define
$$
\aligned
P_{\ge r}(\osX\asD)
&:= H^{r+1}(W^{r}\cK) = \Pic({}^{o}W^{r}X_{\ssb}[r]),
\\
P_{r}(\osX\asD)
&:= H^{r+1}({\Gr}_{W}^{r}\cK) = \Pic(\oX_{r})
\oplus \Gamma (\tD_{r+1},\bZ),
\\
P_{\ge r}(\osX)
&:= H^{r+1}(W^{r}\cK') =
\Pic({}^{o}W^{r}\osX[r]),
\\
P_{r}(\osX)
&:= H^{r+1}({\Gr}_{W}^{r}\cK') =
\Pic(\oX_{r}).
\endaligned
$$
Then we have an exact sequence
$$
0 \to P_{\ge r}(\osX) \to
P_{\ge r}(\osX\asD)
\to \Gamma (\tD_{r+1},\bZ).
$$

By (3.2) below,
$ P_{\ge r}(\osX) $ has a structure of algebraic group
$ \cP_{\ge r}(\osX) $\break (locally of finite type)
such that the identity component is a semiabelian\break variety.
(This is well-known for
$ P_{r}(\osX) $.)
Let
$ P_{\ge r}(\osX)^{0}_{\phantom{|}}$ denote the
identity component of
$ P_{\ge r}(\osX) $.
This is identified with
$ P_{\ge r}(\osX\asD)^{0} $
by the above exact sequence (and similarly for
$ P_{r}(\osX^{0}) $,
$ P_{r}(\osX\asD)^{0} $).

By the boundary map
$ \partial $ of the long exact sequence associated with
$$
\to W^{r-1}\cK \to W^{r-2}\cK
\to {\Gr}_{W}^{r-2}\cK \to
$$
we get a commutative diagram
$$
\matrix
P_{r-2}(\osX)^{0} &\lrar& 
P_{r-2}(\osX\asD) &\lrar& 
P_{r-2}(\osX\asD)/
P_{r-2}(\osX)^{0}
\\ \navs
\bdar{\partial}&&\bdar{\partial}&&\bdar{\partial} 
\\ \navs
P_{\ge r-1}(\osX)^{0} &\lrar& 
P_{\ge r-1}(\osX\asD) &\lrar& 
P_{\ge r-1}(\osX\asD)/
P_{\ge r-1}(\osX)^{0}.
\endmatrix
\leqno(3.1.4)
$$
Let
$ \Gamma'_{r}(\osX\asD) $ be
the kernel of the right vertical morphism of (3.1.4).
Put
$$
\NS(\osX\asD)^{j} :=
P_{j}(\osX\asD
)/P_{j}(\osX)^{0} =\NS(\oX_{j}) \oplus \Gamma
(\tD_{j+1},\bZ),
$$
where
$$
\NS(\oX_{j}) := \Pic(\oX_{j})/\Pic(\oX_{j})^{0}
=\Hom_{\MHS}(\bZ,H^{2}(\oX_{j},\bZ)(1)).
$$
Then
$ (\NS(\osX\asD)^{\ssb},
d^{*}) $ is the single complex associated with a double complex such that
one of the differentials is the Gysin morphism
$$
\Gamma (\tD_{j+1},\bZ) \to \NS(\oX_{j+1}).
$$
Since
$ d^{2} = 0 $,
$ d^{*} $ induces a morphism
$$
d^{*} : \NS(\osX\asD)^{r-3}
\to
\Gamma'_{r}(\osX\asD).
\leqno(3.1.5)
$$

We define
$$
\aligned
\Gamma_{r}(\osX\asD)
&=\Coker(d^{*} : \NS(\osX\asD)^{r-3} \to
\Gamma'_{r}(\osX\asD)),
\\
G_{r}(\osX\asD)
&=\Coker(\partial : \cP_{r-2}(\osX)^{0}
\to \cP_{\ge r-1}(\osX)^{0}),
\endaligned
$$
where
$ \Gamma_{r}(\osX\asD),
\Gamma'_{r}(\osX\asD) $ are
identified with locally finite commutative group schemes.
Then (3.1.4) induces morphisms
$$
\Gamma'_{r}(\osX\asD)
\to G_{r}(\osX\asD
),\quad \Gamma_{r}(\osX\asD)
\to G_{r}(\osX\asD),
$$
which define respectively
$$
M'_{r}(\osX\asD),\quad
M_{r}(\osX\asD).
$$
(This construction is equivalent to the one in [26].)

\vglue12pt {\it Remark.} By (3.1.3),
$ P_{\ge r}(\osX) $ is identified with the group of
isomorphism classes of
$ (L,\gamma) $ where
$ L $ is a line bundle on
$ \oX_{r} $ and
$$
\gamma : \cO_{\oX_{r+1}} \simto d^{*}L
$$
is a trivialization such that
$$
d^{*}\gamma : d^{*}\cO_{\oX_{r+1}}
\,(= \cO_{\oX_{r+2}})
\simto (d^{2})^{*}L = \cO_{\oX_{r+2}}
$$
is the identity morphism.
(Note that
$ d^{*}L $ is defined by using tensor of line bundles.)
See also [2], [26].
\vglue9pt

For the construction of the group scheme
$ \cP_{\ge r}(\osX) $, we need Grothendieck's
theory of representable group functors (see [15], [23]) as follows:

\nonumproclaim{{3.2.~Theorem}}
There exists a
$ k $\/{\rm -}\/group scheme locally of finite type\break
$ \cP_{\ge r}(\osX) $ such that the group of its
$ k $\/{\rm -}\/valued points is isomorphic to
$ P_{\ge r}(\osX) $.
Moreover{\rm ,} \pagebreak
$ \cP_{\ge r}(\osX) $ has the following universal
property\/{\rm :} For any
$ k $\/{\rm -}\/variety
$ S $ and any
$ (L,\gamma) \in P_{\ge r}(\osX\times_{k}S) $ as
above{\rm ,} the set\/{\rm -}\/theoretic map
$ S(k) \to P_{\ge r}(\osX) $ obtained by restricting
$ (L,\gamma) $ to the fiber at
$ s \in S(k) $ comes from a morphism of
$ k $\/{\rm -}\/schemes
$ S \to \cP_{\ge r}(\osX) $.
\endproclaim

\demo{{P}roof}
Essentially the same as in [2].
(This also follows from [24], see Remark after (3.3).)
\enddemo

{\it Remark.} We can easily construct a
$ k $-scheme locally of finite type
$ S $ and
$ (L,\gamma) \in P_{\ge r}(\osX\times_{k}S) $ such
that the associated map
$ f : S(k) \to P_{\ge r}(\osX) $ is surjective by using the theory
of Hilbert scheme.
Then
$ P_{\ge r}(\osX) $ has at most unique structure of
$ k $-algebraic group such that
$ f $ is algebraic.
But it is nontrivial that this really gives an algebraic structure on
$ P_{\ge r}(\osX) $, because it is even unclear if the
inverse image of a closed point is a closed variety for example.
The independence of the choice of
$ (L,\gamma) $ and
$ S $ is also nontrivial.
So we have to use Grothendieck's general theory using sheafification in
the {\it fppf} (faithfully flat and of finite presentation) topology.

If
$ k = \bC $, we can prove (3.2) (for smooth varieties
$ S) $ by using Hodge theory.
See (5.3).
In fact, this implies an isomorphism between the semiabelian parts of
$1$-motives (see also Remark after (5.3)).
But this proof of (3.2) is not algebraic, and cannot be used to prove
Deligne's conjecture.

\nonumproclaim{{3.3.~Lemma}}
The identity component
$ \cP_{\ge r}(\osX)^{0} $ is a semiabelian
variety.
\endproclaim

\demo{{P}roof} With the notation of (3.1), let
$$
T_{r}(\osX) = H^{r}(\Gamma (\osX,
\cO_{{\oX}_{\ssb}}^{*})),\quad
\oP_{r}(\osX) = \Ker(d^{*} :
\Pic(\oX_{r}) \to \Pic(\oX_{r+1})).
$$
Then (3.1.3) induces an exact sequence
$$
0 \to T_{r+1}(\osX) \to P_{\ge
r}(\osX) \to \oP_{r}(\oX
_{\ssb}) \to T_{r+2}(\osX) \to .
\leqno(3.3.1)
$$
Let
$ T_{r}(\osX)^{0}, \oP_{r}(\oX
_{\ssb})^{0} $ denote the identity components of
$ T_{r}(\osX), \oP_{r}
(\osX).
\quad $ Then (3.3.1) induces a short exact sequence
$$
0 \to T_{r+1}(\osX)' \to P_{\ge
r}(\osX)^{0} \to \oP
_{r}(\osX)^{0} \to 0,
\leqno(3.3.2)
$$
where
$ T_{r+1}(\osX)' $ is a subgroup of
$ T_{r+1}(\osX) $ with finite index.
This gives a structure of semiabelian variety
$$
0 \to T_{r+1}(\osX)^{0} \to P_{\ge
r}(\osX)^{0} \to \oP
_{r}(\osX)' \to 0,
\leqno(3.3.3)
$$
with an isogeny of abelian varieties
$$
0 \to T_{r+1}(\osX)'/T_{r+1}(\oX
_{\ssb})^{0} \to \oP_{r}(\osX)'
\to \oP_{r}(\osX)^{0} \to 0.
\leqno(3.3.4)
$$
\enddemo

{\it  Remark}. The representability of the Picard functor follows also
from [24, Prop.\ 17.4], using (3.3.3).
See also [26].

\nonumproclaim{{3.4.~Theorem}}
With the notation of {\rm (3.1),} let
$ W_{2}{H}_{(1)}^{r}(X,Y;\bZ_{l}) $ be the\break
$ \bZ_{l} $\/{\rm -}\/submodule of
$ W_{2}{H}_{\et}^{r}(X,Y;\bZ_{l}) $ whose image in
$ {\Gr}_{W}^{2}{H}_{\et}^{r}(X,Y;\bZ_{l}) $ is generated
by the image of
$ \Gamma'_{r}(\osX^{\phantom{|}} \asD) $ under the cycle map.
Then there is a canonical isomorphism
$$
r_{l}(M_{r}(X,Y))_{\fr}(-1) = W_{2}{H}_{(1)}^{r}(X,Y;\bZ_{l})_{\fr}
$$
compatible with the weight filtration
$ W $.
\endproclaim

\demo{{P}roof}
Note first that
$ W $ is defined by using a resolution
$ \osX $ in (2.3), and
$ X_{\ssb} := \osX \setminus D_{\ssb} $ will be denoted sometimes
by
$ \osX\asD $ as in (3.1).
Let
$$
\cK_{\et} = \Gamma(X_{\ssb},
{C}_{\et}^{\ssb}(\bG_{m})),
$$
where
$ {C}_{\et}^{\ssb} $ denotes the canonical flasque resolution
of Godement in the \'etale topology.
Then
$ \cK_{\et} $ has a filtration
$ W $ in a generalized sense, which is induced by the cofiltration
$ {}^{o}W $ in (3.1) so that
$$
W^{j}\cK_{\et} = \Gamma({X}_{\ssb}^{j},
{C}_{\et}^{\ssb}(\bG_{m})),
$$
where
$ {X}_{\ssb}^{j} := {}^{o}W^{j}\osX\asD $.
We define
$$
\cK_{\et}(r,n) =
C(W^{r-1}\cK_{\et} \overset{n}\to\to
W^{r-2}\cK_{\et}),
$$
and similarly for
$ \cK(r,n) $ in the Zariski topology.

By the Kummer sequence
$ 0 \to \mu_{n} \to \bG_{m} \overset{n}\to\to
\bG_{m} \to 0, $ we have
$$
{H}_{\et}^{i}({X}_{\ssb}^{j},\mu_{n}) =
H^{i-1}(C(W^{j}\cK_{\et} \overset{n}\to\to W^{j}\cK_{\et})).
$$
Note that the \'etale cohomology
$ {H}_{\et}^{i}({X}_{\ssb}^{j},\mu_{n}) $ has
the weight filtration
$ W $, and
$$
W_{q-2}{H}_{\et}^{i}({X}_{\ssb}^{j},\mu_{n}) =
\Im({H}_{\et}^{i}({X}_{\ssb}^{i-q},\mu_{n}) \to
{H}_{\et}^{i}({X}_{\ssb}^{j},\mu_{n}))
\leqno(3.4.1)
$$
for
$ q \le 2 $ and
$ i - q \ge j $ as in (4.4.2).
Here the shift of
$ W $ by
$ 2 $ comes from the Tate twist
$ \mu_{n} $.
We define
$$
\aligned
E'(r,n)
&= \Im(H^{r-1}\cK_{\et}(r,n) \to
{H}_{\et}^{r}
({X}_{\ssb}^{r-2},\mu_{n})),
\\
N(r,n)
&= \Im(\partial : H^{r-1}{\Gr}_{W}^{r-2}\cK_{\et}
\to {H}_{\et}^{r}(X_{\ssb}^{r-2}/X_{\ssb}^{r-1},\mu_{n}))
\\
(&= \Coker(n : H^{r-1}{\Gr}_{W}^{r-2}\cK_{\et}
\to H^{r-1}{\Gr}_{W}^{r-2}\cK_{\et})),
\endaligned
$$
where
$ {\Gr}_{W}^{r-2}\cK_{\et} $ and
$ X_{\ssb}^{r-2}/X_{\ssb}^{r-1} := {\Gr}_{{}^{o}W}^{r-2}
\osX\asD $ are defined by using the mapping cones, and
$ \partial $ is induced by the boundary map of the Kummer sequence, and
gives the cycle map.
Note that
$$
\aligned
H^{r-1}{\Gr}_{W}^{r-2}\cK_{\et}
&=\Pic(\oX_{r-2}) \oplus
\Gamma(\tD_{r-1},\bZ),
\\
{H}_{\et}^{r}(X_{\ssb}^{r-2}/X_{\ssb}^{r-1},\mu_{n})
&={H}_{\et}^{2}(\oX_{r-2},\mu_{n}) \oplus
H^{0}(\tD_{r-1},\bZ/n),
\endaligned
$$
using Hilbert's theorem 90 for the first and the local cohomology for
the second.
Let
$$
N'(r,n) = N(r,n) \cap \ker({H}_{\et}^{r}
(X_{\ssb}^{r-2}/X_{\ssb}^{r-1},\mu_{n})
\to {H}_{\et}^{r+1}({X}_{\ssb}^{r-1},\mu_{n})).
$$
This coincides with the intersection of
$ N(r,n) $ with the image of
$ {H}_{\et}^{r}({X}_{\ssb}^{r-2},\mu_{n}) $ using the long exact
sequence.
So we get a short exact sequence
$$
0 \to W_{-1}{H}_{\et}^{r}
({X}_{\ssb}^{r-2},\mu_{n}) \to E'(r,n) \to
N'(r,n) \to 0
\leqno(3.4.2)
$$
by considering the natural morphism between the distinguished
triangles
$$
\matrix
 \to &  C(n : W^{r-1}\cK_{\et})
&\to &  \cK_{\et}(r,n)
&\to &  {\Gr}_{W}^{r-2}\cK_{\et} &\to & 
\\
 &\Big\Vert&&\Big\downarrow&&\Big\downarrow
\\
 \to &  C(n : W^{r-1}\cK_{\et})
&\to &  C(n : W^{r-2}\cK_{\et})
&\to &  C(n : {\Gr}_{W}^{r-2}\cK_{\et}) &\to & ,
\endmatrix
$$
where
$ C(n : A) $ is the abbreviation of
$ C(n : A \to A) $ for an abelian group
$ A $.

Let
$ E'(r,l^{\infty}) $ be the projective limit of
$ E'(r,l^{m}) $.
By the Mittag-Leffler condition, it is identified with a
$ \bZ_{l} $-submodule of
$ {H}_{\et}^{r}({X}_{\ssb}^{r-2},\bZ_{l}) $ so
that
$$
W_{-1}{H}_{\et}^{r}({X}_{\ssb}^{r-2},\bZ_{l})
\subset E'(r,l^{\infty}) \subset
{H}_{\et}^{r}({X}_{\ssb}^{r-2},\bZ_{l})
$$
and
$ E'(r,l^{\infty})/W_{-1}{H}_{\et}^{r}({X}_{\ssb}^{r-2},\bZ_{l})
\subset{\Gr}_{W}^{0}{H}_{\et}^{r}({X}_{\ssb}^{r-2},\bZ_{l}) $ is
generated by the images of
$ \Gamma'_{r}(\osX\asD) $ under the cycle map.
Let
$$
E(r,l^{\infty}) = \Im(E'(r,l^{\infty}) \to
H_{\et}^{r}(X_{\ssb}^{r-3},\bZ_{l})).
$$
Since the integral weight spectral sequence degenerates at
$ E_{2} $ modulo torsion, we get
$$
W_{2}{H}_{(1)}^{r}(X,Y;\bZ_{l})_{\fr} = E(r,l^{\infty})_{\fr}.
\leqno(3.4.3)
$$

We have to show that
$ E(r,l^{\infty})_{\fr} $ is naturally isomorphic to the
$ l $-adic realization of
$ M_{r}(\osX\asD)_{\fr} $.
Define a decreasing filtration
$ G $ on
$$C(W^{r-2}\cK_{\et} \overset{n}\to\to W^{r-2}
\cK_{\et})$$ by
$$
\alignat 2
G^{0}
&= C(W^{r-2}\cK_{\et} \overset{n}\to\to W^{r-2}\cK_{\et}),
&\quad G^{1}
&=\cK_{\et}(r,n),
\\
G^{2}
&= C(0 \to W^{r-1}\cK_{\et}),
&\quad G^{3}
&= 0.
\endalignat
$$
Then
$ {\Gr}_{G}^{0} = {\Gr}_{W}^{r-2}\cK_{\et}[1], {\Gr}_{G}^{1}
= W^{r-1}\cK_{\et}[1]
\oplus {\Gr}_{W}^{r-2}\cK_{\et} $.
So there is a canonical morphism in the derived category
$ \beta : {\Gr}_{W}^{r-2}\cK_{\et} \to
\cK_{\et}(r,n) $ whose mapping cone is isomorphic to
$ C(W^{r-2}\cK_{\et} \overset{n}\to\to W^{r-2}
\cK_{\et}) $.
By the associated long exact sequence, we have
$$
E'(r,n) = \Coker(\beta : H^{r-1}{\Gr}_{W}^{r}\cK_{\et}
\to H^{r-1}\cK_{\et}(r,n)).
\leqno(3.4.4)
$$
Here we can replace
$ \cK_{\et} $ with
$ \cK $,
because the right-hand side does not change by doing it.

Using the induced filtration
$ G $ on
$ \cK(r,n) $,
we have a long exact sequence
$$
\aligned
H^{r-1}(W^{r-1}\cK)
&\oplus H^{r-2}({\Gr}_{W}^{r-2}\cK)
\to H^{r-1}(W^{r-1}\cK)
\\
\to H^{r-1}\cK(r,n) \to H^{r}(W^{r-1}\cK)
&\oplus H^{r-1}({\Gr}_{W}^{r-2}\cK) \to
H^{r}(W^{r-1}\cK),
\endaligned
$$
where the first and the last morphisms are given by the sum of the
multiplication by
$ n $ and the boundary morphism
$ \partial $.
In particular, the cokernel of the first morphism is a finite group, and
is independent of
$ n = l^{m} $ for
$ m $ sufficiently large, because
$$
H^{r-1}(W^{r-1}\cK) = \Ker(d^{*} :
\Gamma(\oX_{r-1},\cO_{\oX_{r-1}}^{*}) \to
\Gamma(\oX_{r},\cO_{\oX_{r}}^{*})).
$$

Consider the cohomology
$ H(r,n) $ of
$$
H^{r-1}({\Gr}_{W}^{r-2}\cK) \to H^{r}(W^{r-1}\cK) \oplus
H^{r-1}({\Gr}_{W}^{r-2}\cK) \to H^{r}(W^{r-1}\cK),
$$
where the first morphism is the composition of
$ \beta $ with the third morphism of the above long exact sequence.
Then, by the above argument, there is a surjective canonical morphism
$$
E'(r,n) \to H(r,n),
\leqno(3.4.5)
$$
whose kernel is a finite group, and is independent of
$ n = l^{m} $ for
$ m $ sufficiently large.
By definition,
$ H(r,n) $ is isomorphic to
$ H^{0} $ of
$$
C(\partial : H^{r-1}({\Gr}_{W}^{r-2}\cK) \to
H^{r}(W^{r-1}\cK))\otimes C(n : \bZ \to \bZ)[-1].
$$
Since the Tate module of a finitely generated abelian group vanishes,
we can replace the mapping cone of
$ \partial : H^{r-1}({\Gr}_{W}^{r-2}\cK) \to
H^{r}(W^{r-1}\cK) $ with that of
$$
\partial : \Ker(H^{r-1}({\Gr}_{W}^{r-2}\cK) \to
P_{\ge r-1}(\osX)/P_{\ge r-1}
(\osX)^{0}) \to
P_{\ge r-1}(\osX)^{0}
$$
in the notation of (3.1).
Furthermore,
$ H^{0} $ of
$$
C(\partial : P_{r-2}(\osX)^{0} \to
\partial (P_{r-2}(\osX)^{0}))\otimes C(n : \bZ
\to \bZ)[-1]
$$
is a finite group, and is independent of
$ n = l^{m} $ for
$ m $ sufficiently large.

On the other hand, the
$ l $-adic realization
$ r_{l}(M'_{r}(\osX\asD)) $ is
the projective limit of
$ H^{0} $ of
$$
C(\Gamma'_{r}(\osX\asD)
\to G_{r}(\osX\asD))
\otimes C(n : \bZ \to \bZ)[-1].
$$
So we get a canonical surjective morphism
$$
E'(r,l^{\infty}) \to r_{l}(M'_{r}(\osX\asD))
$$
whose kernel is a finite group.
This is clearly compatible with the weight filtration.
Then the assertion follows by taking the image in
$ H_{\et}^{r}(X_{\ssb},\bZ_{l})_{\fr} $,
and using the
$ E_{2} $-degeneration of the weight spectral sequence modulo torsion.
\enddemo

\nonumproclaim{{3.5.~Theorem}}
With the notation of {\rm (3.1),} let
$ {H}_{\DR,(1)}^{r}(X,Y) $ be the\break
$ k $\/{\rm -}\/submodule of
$ W_{2}{H}_{\DR}^{r}(X,Y) $ whose image in
$ {\Gr}_{W}^{2}{H}_{\DR}^{r}(X,Y) $ is generated by the \pagebreak  image of
$ \Gamma'_{r}(\osX^{\phantom{|}}\asD) $ under the cycle map{\rm ,} where
$ {H}_{\DR}^{r}(X,Y) $ is defined by using
$ \osX\asD $ in {\rm (2.3).}
Then there is a canonical isomorphism
$$
r_{\DR}(M_{r}(X,Y))(-1) = {H}_{\DR,(1)}^{r}(X,Y).
$$
compatible with the Hodge filtration
$ F $ and the weight filtration
$ W $.
\endproclaim

\demo{{P}roof}
By definition of de Rham realization [10], we have to construct first the
universal
$ \bG_{a} $-extension of
$ M_{r}(X,Y) = M_{r}(\osX\asD) $.
Let
$ \tcK $ be the shifted mapping cone
$ [\cK \to \cK^{1}] $ with
$$
\cK = \Gamma (\osX,
\cC^{\ssb}(j_{*}\cO_{{X}_{\ssb}}^{*})),\quad
\cK^{1} = \Gamma (\osX,
\cC^{\ssb}({\Omega}_{{\oX}_{\ssb}}^{1}(\log D_{\ssb}))),
$$
where
$ \cK $ has  degree zero, and the morphism is induced by
$ g \mapsto g^{-1}dg $.
Here
$ \cC^{\ssb} $ is the canonical flasque resolution of Godement in the
Zariski
topology.
Then
$ \tcK $ has a filtration
$ W $ in a generalized sense, which is induced by the cofiltration
$ {}^{o}W $ on
$ \osX\asD $.

Let
$ G $ be the convolution of
$ W $ with the Hodge filtration
$ F $ defined by
$ {\Gr}_{F}^{0} = \cK $ and
$ {\Gr}_{F}^{1} = \cK^{1} $.
Then
$ G^{i}\tcK = [W^{i}\cK \to
W^{i-1}\cK^{1}] $ with
$$
W^{i}\cK = \Gamma ({\oX}_{\ssb}^{i},
\cC^{\ssb}({j}_{*}^{i}\cO_{{X}_{\ssb}^{i}}^{*})),\quad
W^{i-1}\cK^{1} = \Gamma ({\oX}_{\ssb}^{i-1},
\cC^{\ssb}({\Omega}_{{\oX}_{\ssb}^{i-1}}^{1}
(\log {D}_{\ssb}^{i-1}))),
$$
where
$ {X}_{\ssb}^{i} = {}^{o}W^{i}\osX\asD, {D}_{\ssb}^{i} =
{\oX}_{\ssb}^{i} \setminus {X}_{\ssb}^{i} $,
and
$ {\oX}_{\ssb}^{i} $ is the closure of
$ {X}_{\ssb}^{i} $ in
$ {\oX}_{\ssb} $ with the inclusion
$ j^{i} : {X}_{\ssb}^{i} \to {\oX}_{\ssb}^{i} $.
This is compatible with
$ W $ on
$ \cK $ in (3.1).
Note that
$ {\Omega}_{{\oX}_{j}^{i}}^{1}(\log {D}_{j}^{i}) =
{\Omega}_{{\oX}_{j}}^{1}(\log D_{j}) $ for
$ j > i $,
$ {\Omega}_{{\oX}_{j}}^{1} $ for
$ j = i $,
and
$ 0 $ for
$ j < i $.
Let
$$
\cK' = \Gamma (\osX, \cC
(\cO_{\osX})),\quad W^{i}\cK' =
\Gamma ({\oX}_{\ssb}^{i},
\cC^{\ssb}(\cO_{{\oX}_{\ssb}^{i}})),
$$
and define
$ W, F, G $ similarly on
$ \tcK' $ so that
$ G^{i}\tcK' = [W^{i}\cK' \to
W^{i-1}\cK^{1}] $, where
$ \cK^{1} $ and
$ W^{i-1}\cK^{1} $ are as above, and the morphism is induced by
$ d $.

Then we have an exact sequence
$$
0 \to H^{r-1}W^{r-2}\cK^{1} \to
H^{r}G^{r-1}\tcK \overset{\partial}\to\to
H^{r}W^{r-1}\cK \to,
\leqno(3.5.1)
$$
and
$ \Im\,\partial $ contains
$ P_{\ge r-1}(\osX)^{0} $ with the notation of (3.1).
Indeed, this is reduced to the case
$ k = \bC $,
and is verified by using the canonical morphism of the corresponding
exact sequence
$$
0 \to H^{r-1}W^{r-2}\cK^{1} \to
H^{r}G^{r-1}\tcK' \overset{\partial}\to\to
H^{r}W^{r-1}\cK' \to 0
$$
to the above sequence, because
$ d : H^{j}W^{i}\cK' \to H^{j}W^{i-1}
\cK^{1} $ vanishes by Hodge theory and the torsion of
$ H^{r-1}W^{r-1}\cK $ comes from
$ H^{r-1}G^{r-1}\tcK $.

Let
$ (H^{r}G^{r-1}\tcK)^{0} $ be the
$ k $-submodule of
$ H^{r}G^{r-1}\tcK $ which contains\break
$ H^{r-1}W^{r-2}\cK^{1} $,
and whose image by
$ \partial $ is
$ P_{\ge r-1}(\osX)^{0} $.
In the case
$ k = \bC $,
this is the image of
$ H^{r}G^{r-1}\tcK' $.
We can verify that
$ (H^{r}G^{r-1}\tcK)^{0} $ has naturally a structure of a commutative
$ k $-group scheme.
We consider the image of
$ (H^{r}G^{r-1}\tcK)^{0} $ by the canonical morphism
$$
H^{r}G^{r-1}\tcK \to H^{r-1}{\Gr}_{W}^{r-2}\cK^{1}.
$$
This coincides with the image of
$ H^{r-1}W^{r-2}\cK^{1} $.
Indeed, we can replace\break
$ H^{r}G^{r-1}\tcK $ with
$ H^{r}G^{r-1}\tcK' $ by reducing to the case
$ k = \bC $.

Let
$ (H^{r-1}{\Gr}_{W}^{r-2}\cK^{1})_{\alg} $ be the
$ k $-submodule of
$$
H^{r-1}{\Gr}_{W}^{r-2}\cK^{1} = H^{1}(\oX_{r-2},
{\Omega}_{{\oX}_{r-2}}^{1})\oplus
\Gamma (\tD_{r-1},\cO_{\tD_{r-1}})
$$
generated by the divisor classes
in
$ H^{1}(\oX_{r-2}, {\Omega}_{{\oX}_{r-2}}^{1}) $ and by
$ \Gamma (\tD_{r-1},\cO_{\tD_{r-1}}) $.
Let
$ (H^{r}G^{r-1}\tcK{)}_{\alg}^{0} $ be the largest
$ k $-submodule of
$ (H^{r}G^{r-1}\tcK)^{0} $ whose image in
$ H^{r-1}{\Gr}_{W}^{r-2}\cK^{1} $ is contained in
$ (H^{r-1}{\Gr}_{W}^{r-2}\cK^{1})_{\alg} $.
We define similarly
$$
(H^{r}G^{r-1}\tcK')_{\alg} \subset H^{r}G^{r-1}\tcK',\quad
(H^{r-1}W^{r-2}\cK^{1})_{\alg} \subset H^{r-1}W^{r-2}
\cK^{1}.
$$
These have the induced filtrations
$ F $ and
$ W $ so that
$$
\aligned
{\Gr}_{F}^{0}(H^{r}G^{r-1}\tcK{)}_{\alg}^{0} =
P_{\ge r-1}(\osX)^{0},\quad
&{\Gr}_{F}^{1}(H^{r}G^{r-1}\tcK{)}_{\alg}^{0} =
(H^{r-1}W^{r-2}\cK^{1})_{\alg},
\\
{\Gr}_{F}^{0}(H^{r}G^{r-1}\tcK')_{\alg} =
H^{r}W^{r-1}\cK',\quad
&{\Gr}_{F}^{1}(H^{r}G^{r-1}\tcK')_{\alg} =
(H^{r-1}W^{r-2}\cK^{1})_{\alg},
\endaligned
$$
and
$ W $ on these spaces is calculated by using the weight spectral
sequence which degenerates at
$ E_{2} $,
see (4.4.2).

Consider the morphism
$$
H^{r-1}{\Gr}_{W}^{r-2}\cK \to H^{r}G^{r-1}\tcK
$$
induced by the distinguished triangle
$$
\to G^{r-1}\tcK \to [W^{r-2}\cK
\to W^{r-2}\cK^{1}] \to {\Gr}_{W}^{r-2}\cK \to.
$$
Let
$ (H^{r-1}{\Gr}_{W}^{r-2}\cK)^{(0)} $ be the kernel of the morphism
to
$ H^{r}W^{r-1}\cK/P_{\ge r-1}(\osX)^{0} $.
Then the above morphism induces
$$
(H^{r-1}{\Gr}_{W}^{r-2}\cK)^{(0)} \to
(H^{r}G^{r-1}\tcK)_{\alg}^{0}.
\leqno(3.5.2)
$$
We divide the source by
$ \Pic(\oX_{r-2})^{0}_{\phantom{|}} $ and the target by its image.
(Note\break that the image of
$ \Pic(\oX_{r-2})^{0}$ in
$ (H^{r}G^{r-1}\tcK)_{\alg}^{0} $ is isomorphic to that in\break
$ \Gr_{F}^{0}(H^{r}G^{r-1}\tcK)_{\alg}^{0} $ because there is no
nontrivial morphism of an abelian\break variety to an affine space.)
Since the source is then
$ \Gamma'_{r}(\osX\asD) $,
we can divide these further by the images of
$ \NS(\osX\asD)^{r-3} $ and
$ (H^{r-2}\Gr_{W}^{r-3}\cK^{1})_{\alg} $
in the notation of (3.1).
Let
$ \tG_{r}(\osX\asD) $ be the commutative
$ k $-group scheme whose underlying group of
$ k $-valued points is the cokernel of the canonical morphism
$$
\Pic(\oX_{r-2})^{0} \oplus (H^{r-2}\Gr_{W}^{r-3}\cK^{1})_{\alg}
\to (H^{r}G^{r-1}\tcK)_{\alg}^{0},
$$
which underlies naturally a morphism of groups schemes.
Then we get
$$
\tM_{r}(\osX\asD) :=
[\Gamma_{r}(\osX\asD)\to \tG_{r}(\osX\asD)].
$$
The Lie algebra
$ \Lie \tG_{r}(\osX\asD) $ is isomorphic to
$$
\Coker(H^{1}(\oX_{r-2},\cO_{{\oX}_{r-2}}) \oplus
(H^{r-2}\Gr_{W}^{r-3}\cK^{1})_{\alg}
\to (H^{r}G^{r-1}\tcK')_{\alg})
$$
by the standard argument using
$ \Spec\, k[\varepsilon ] $ (cf. [21]).
But this is isomorphic to the image of
$ (H^{r}G^{r-1}\tcK')_{\alg} $ in
$ {H}_{\DR}^{r}(X_{\ssb})/F^{2} $
by the
$ E_{2} $-degeneration of the weight spectral sequence together with the
strictness of the Hodge filtration.\break\vglue-13pt\noindent 
Since it is isomorphic to
$ {H}_{\DR,(1)}^{r}(X_{\ssb}) $, it is sufficient
to show that
$ \tM_{r}(\osX\asD) $ is
the universal
$ \bG_{a} $-extension of
$ M_{r}(\osX\asD) $.
Then we may replace it with
$ {\Gr}_{W}^{i} $,
and the assertion is reduced to the well-known fact about the universal
$ \bG_{a} $-extension of the Picard variety (see Remark (i) below).
This completes the proof of   (3.5).
\enddemo

{\it Remarks.} (i)
Let
$ \tM = [\Gamma \to \tG] $ be the universal
$ \bG_{a} $-extension of a\break $1$-motive
$ M = [\Gamma \to G] $.
Then the de Rham realization
$ r_{\DR}(M) $ is defined to be
$ \Lie \tG $.
It is known that the universal extension is given by a commutative
diagram with exact rows
$$
\matrix
&&&&\Gamma &\dlines& \Gamma
\\ \navs
&&&&\Big\downarrow&&\Big\downarrow
\\ \navs
0 &\lrar&  \Ext^{1}(M,\bG_{a})^{\vee} &\lrar&  \tG &\lrar& 
G &\lrar&  0
\endmatrix
$$
where
$ \Ext^{1}(M,\bG_{a}) $ is a finite dimensional
$ k $-vector space, and its dual is identified with a group scheme.
(See [10], [21] and also [2].)
Indeed, if
$ 0 \to V \to M' \to M \to 0 $ is an extension by a
$ k $-vector space
$ V $ which is identified with a
$ k $-group scheme, it gives a morphism
$$V^{\vee}\,(= \Hom(V,\bG_{a})) \to \Ext^{1}(M,\bG_{a})$$ by
composition, and its dual deduces the original
extension from the universal extension.
In particular, the functor
$ M \to \tM $ is exact.
If
$ M $ is a torus,
$ \tM = M $.
If
$ M = \Pic(X)^{0} $ for a smooth proper variety, then
$ \Lie \tG = {H}_{\DR}^{1}(X) $ and
$ \Ext^{1}(M,\bG_{a})^{\vee} = F^{1}{H}_{\DR}^{1}(X) = \Gamma
(X,{\Omega}_{X}^{1}) $.
If
$ M = [\Gamma \to 0] $,
then
$ \tM = [\Gamma \to \Gamma \otimes
\bG_{a}] $.
(See {\it loc.\ cit}.)
\vglue5pt

(ii) With the above notation, assume
$ k = \bC $.
Then we have a commutative diagram with exact rows
$$
\matrix
0 &\lrar&  E &\lrar&  E' &\lrar&  \Gamma_{r}
(\osX\asD) &\lrar&  0
\\ \navs
&&\Big\Vert&&\Big\downarrow&&\Big\downarrow
\\ \navs
0 &\lrar&  E &\lrar&  \Lie G_{r}(\osX\asD)& \raise3pt\hbox{${{\exp}\atop{\textstyle\lrar}}$}&
G_{r}(\osX\asD)
&\lrar&  0
\endmatrix
$$
where
$ \exp $ is the exponential map.
Note that the exponential map of a commutative Lie group depends on the
analytic structure of the group, and it cannot be used for the proof of
the coincidence of the natural analytic structure with the one coming
from the algebraic structure of the Picard variety,
before we show that the generalized Abel-Jacobi map
depends analytically on the parameter.
See also Remark  \pagebreak after (5.3). 

We will later show that
$ E' $ is identified with
$ {H}_{(1)}^{r}(X,Y;\bZ) $ modulo torsion.
This is true if and only if we have the above commutative diagram with
$ E' $ replaced by
$ {H}_{(1)}^{r}(X,Y;\bZ) $ (modulo torsion).
But it is easy to see that the assertion is equivalent to the coincidence
of the two extension classes associated with the $1$-motive
$ M_{r}(X,Y) $ and the mixed Hodge structure
$ {H}_{(1)}^{r}(X,Y;\bZ) $.
This will be proved in the proof of (5.4).
This point is not clear   in [26].

\section{Mixed Hodge theory}

We review the theory of mixed Hodge complexes ([10], [4]) and Deligne
cohomology ([3], [4]).
See also [11], [12], [13], [19], etc.
We assume
$ k = \bC $.

\vglue12pt  4.1. {\it Mixed Hodge complexes.} Let
$ \cC_{\cH} $ be the category of mixed Hodge complexes in the
sense of Beilinson [4, 3.2].
An object
$ K \in \cC_{\cH} $ consists of (filtered or bifiltered)
complexes
$ K_{\bZ}, K'_{\bQ}, (K_{\bQ},W), (K'_{\bC},W), (K
_{\bC};F,W) $ over
$ \bZ, \bQ $ or
$ \bC $, together with (filtered) morphisms
$$
\alignat 2
&\alpha_{1} : K_{\bZ} \to K'_{\bQ},\quad
&&\alpha_{2} : K_{\bQ} \to K'_{\bQ},
\\
&\alpha_{3} : (K_{\bQ}, W) \to (K'_{\bC},W),\quad
&&\alpha_{4} : (K_{\bC},W) \to (K'_{\bC},W),
\endalignat
$$
which induce quasi-isomorphisms after scalar extensions.
These complexes are bounded below, the
$ H^{j}K_{\bZ} $ are finite
$ \bZ $-modules and vanish for
$ j \gg 0 $, the filtration
$ F $ on
$ {\Gr}_{i}^{W}K_{\bC} $ is strict, and
$ H^{j}{\Gr}_{i}^{W}(K_{\bQ},(K_{\bC},F)) $ is a pure Hodge
structure of weight
$ i $,
using the isomorphism
$ H^{j}{\Gr}_{i}^{W}K_{\bQ}\otimes_{\bQ}\bC =
H^{j}{\Gr}_{i}^{W} K_{\bC} $ given by
$ \alpha_{3} $ and
$ \alpha_{4} $.
A morphism of
$ \cC_{\cH} $ is a family of morphisms of (filtered or
bifiltered) complexes compatible with the
$ \alpha_{i} $.
A homotopy is defined similarly.
We get
$ \cD_{\cH} $ by inverting bifiltered quasi-isomorphisms.
See {\it loc.~cit.}\ for details.

Similarly, we have categories
$ \cC_{\cH^{p}}, \cD_{\cH^{p}} $ of mixed
$ p $-Hodge complexes.
This is defined by modifying the above definition as follows:

Firstly, the weight of
$ H^{j}{\Gr}_{i}^{W}(K_{\bQ},(K_{\bC},F)) $ is
$ i+j $ (as in [10]) instead of
$ i $, and it is assumed to be polarizable.
A homotopy
$ h $ preserves the Hodge filtration
$ F $.
But it preserves the weight filtration
$ W $ up to the shift
$ -1 $, and
$ dh + hd $ preserves
$ W $.
(This is necessary to show the acyclicity of the mapping cone of the
identity.)
The derived category
$ \cD_{\cH^{p}} $ is obtained by inverting quasi-isomorphisms
(preserving
$ F, W) $.

We have natural functors
$$
\Dec : \cC_{\cH^{p}} \to \cC_{\cH},
\quad \Dec : \cD_{\cH^{p}} \to \cD_{\cH},
\leqno(4.1.1)
$$
by replacing the weight filtration
$ W $ with
$ \Dec\, W $ (see [10]), which is defined by
$$
(\Dec\, W)_{j}K_{\bQ}^{i} =
\Ker(d : W_{j-i}K_{\bQ}^{i}
\to {\Gr}_{j-i}^{W}K_{\bQ}^{i+1}).
$$
Note that (4.1.1) is well-defined, because
$ (F, \Dec\, W) $ is bistrict, and a quasi-isomorphism (preserving
$ F, W) $ induces a bifiltered quasi-isomorphism for\break
$ (F, \Dec\, W) $.
See [27, 1.3.8] or [30, A.2].

The Tate twist
$ K(m) $ of
$ K \in \cD_{\cH} $ (or
$ \cD_{\cH^{p}}) $ for
$ n \in \bZ $ is defined by twisting the complexes over
$ \bZ $ or
$ \bQ $ and shifting the Hodge filtration
$ F $ and the weight filtration
$ W $ as usual [10]; e.g.\
$ W_{j}(K(m)_{\bQ}) = W_{j+2m}K_{\bQ}(m), F^{p}(K(m)_{\bC})
=F^{p+m}K_{\bC} $.

For
$ K \in \cD_{\cH} $ we define
$ \Gamma_{\cD}K $ and
$ \Gamma_{\cH}K $ using the shifted mapping cones (i.e.\  the first
terms have degree zero):
$$
\aligned
\Gamma_{\cD}K
&= [K_{\bZ}\oplus K_{\bQ}\oplus
F^{0}K_{\bC}
\to K'_{\bQ}\oplus K'_{\bC}],
\\
\Gamma_{\cH}K
&= [K_{\bZ}\oplus W_{0}K_{\bQ}\oplus
F^{0}W_{0}K_{\bC}
\to K'_{\bQ}\oplus W_{0}K'_{\bC}],
\endaligned
$$
where the morphisms of complexes are given by
$ (a,b,c) \mapsto (\alpha_{1}(a) - \alpha_{1}(b)$,\break
$\alpha_{1}(b)-\alpha_{1}(c)) $.
Note that we have a quasi-isomorphism
$$
\Gamma_{\cD}K \to [K_{\bZ}\oplus F^{0}K_{\bC}
\to K'_{\bC}],
\leqno(4.1.2)
$$
where the morphism of complexes is given by
$ \alpha'\scirc \alpha_{1} - \alpha_{4} $, if there is a morphism
$$
\alpha' : K'_{\bQ} \to K'_{\bC}
$$
such that
$ \alpha'\scirc \alpha_{2} = \alpha_{3} $.
(Indeed, the quasi-isomorphism is given by
$ (a,b,c;b',c')\break \mapsto (a,c;b'+c') $.)
We can also define a similar
complex for polarizable mixed Hodge complexes.
But it is not used in this paper.
For
$ K \in \cD_{\cH^{p}} $, we define
$$
\Gamma_{\cD}K = \Gamma_{\cD}(\Dec\, K),
\quad \Gamma_{\cH}K =\Gamma_{\cH}(\Dec\, K),
$$
using
$ \Dec $ in (4.1.1).

By Beilinson [4, 3.6], we have a canonical isomorphism
$$
\Hom_{\cH}(\bZ,K) = H^{0}(\Gamma_{\cH}K),
\leqno(4.1.3)
$$
where
$ \Hom_{\cH} $ means the group of morphisms in
$ \cD_{\cH} $.
He also shows ({\it loc.~cit.}\/, 3.4) that the canonical functor induces an
equivalence of categories
$$
D^{b}\MHS \simto \cD_{\cH},
\leqno(4.1.4)
$$
where the source is the bounded derived category of mixed
$ \bZ $-Hodge structures.

\vglue12pt 4.2. {\it Deligne cohomology.} For a smooth complex algebraic
variety
$ X $, let
$ \oX $ be a smooth compactification such that
$ D := \oX \setminus X $ is a divisor with simple normal
crossings.
Let
$ j : X \to \oX $ denote the inclusion, and
$ {\Omega}_{{\oX}^{\an}}^{\ssb}\aD $ the
complex of holomorphic logarithmic forms with the Hodge filtration
$ F $ (defined by
$ \sigma) $ and the weight filtration
$ W $.
See [10].
Let
$ \cC^{\ssb} $ denote the canonical flasque resolution of
Godement.
Then we define the mixed Hodge complex associated with
$ (\oX,D) $:
$$
K_{\cH^{p}}(\oX\aD) \in
\cC_{\cH^{p}}
$$
as follows (see [4], [10]): Let
$$
\aligned
&K_{\bZ} = \Gamma (X^{\an},\cC^{\ssb}(\bZ_{X^{\an}})),\quad
K_{\bQ} = K'_{\bQ} = \Gamma(\oX^{\an},j_{*}\cC^{\ssb}(\bQ_{X^{\an}})),
\\
&K_{\bC} = \Gamma (\oX^{\an},{\Omega}_{{\oX}^{\an}}^{\ssb}\aD) ,
\endaligned
$$
where
$ W $ on
$ K_{\bQ} $ and
$ K_{\bC} $ is induced by
$ \tau $ on
$ j_{*}\cC^{\ssb}(\bQ_{X^{\an}}) $ and
$ W $ on
$ {\Omega}_{{\oX}^{\an}}^{\ssb}\aD $
respectively, and the Hodge filtration
$ F $ is induced by
$ \sigma $ on
$ {\Omega}_{{\oX}^{\an}}^{\ssb}\aD $.
We define
$ (K'_{\bC},W) $ by taking the global section functor of the mapping
cone of
$$
(\cC^{\ssb}({\Omega}_{{\oX}^{\an}}^{\ssb}\aD),\tau)\to
(j_{*}\cC^{\ssb}({\Omega}_{{X}^{\an}}^{\ssb}),\tau)\oplus
(\cC^{\ssb}({\Omega}_{{\oX}^{\an}}^{\ssb}\aD),W).
\leqno(4.2.1)
$$
(The description in [30, 3.3] is not precise.
We need a mapping cone as above.)
Note that
$ K_{\cH^{p}}(\oX\aD) $ has the weight
filtration
$ W $ defined over
$ \bZ $.

We will denote the image of
$ K_{\cH^{p}}(\oX\aD) $ in
$ \cD_{\cH^{p}} $ by
$$
K_{\cH^{p}}(X) \in \cD_{\cH^{p}},
$$
because it is independent of the choice of the compactification
$ \oX $ by definition of
$ \cD_{\cH^{p}} $.
We define
$$
K_{\cH}(X) \in \cD_{\cH}
$$
to be the image of
$ K_{\cH^{p}}(X) $ by (4.1.1).

Let
$ X_{\ssb}, \osX $ be as in (3.1).
Applying the above construction to each
$ \oX_{j}, D_{j} $, we get
$$
K_{\cH^{p}}(\osX\asD)
\in \cC_{\cH^{p}}\quad \text{and\quad} K_{\cH^{p}}
(X_{\ssb}) \in \cD_{\cH^{p}},
$$
Here the filtration
$ W $ for
$ K_{\cH^{p}}(\oX_{j}\la D_{j}\ra) $ is shifted by
$ -j $ when the complex is shifted.
We define
$$
K_{\cH}(X_{\ssb}) \in \cD_{\cH}
$$
to be the image of
$ K_{\cH^{p}}(X_{\ssb}) $ by (4.1.1).
Here
$ W $ is not shifted depending on
$ j $, because we take
$ \Dec $.
Then we have a canonical isomorphism of mixed Hodge structures
$$
H^{i}(X_{\ssb}) = H^{i}K_{\cH}(X_{\ssb}).
\leqno(4.2.2)
$$

We define
$$
\aligned
{H}_{\cD}^{i}(X_{\ssb},\bZ(j))
&= H^{i}\Gamma_{\cD}(K_{\cH}(X_{\ssb})(j)),
\\
{H}_{\AH}^{i}(X_{\ssb},\bZ(j))
&= H^{i}\Gamma_{\cH}(K_{\cH}(X_{\ssb})(j)).
\endaligned
$$

For a closed subvariety
$ Y $ of
$ X $, we apply the above construction to a resolution of
$ [\oY \to \oX] $ as in the proof of (2.3), and
get
$$
K_{\cH^{p}}(X,Y) \in \cD_{\cH^{p}},\quad K_{\cH}(X,Y)
\in \cD_{\cH}.
$$
These are independent of the choice of the resolution by definition of
$ \cD_{\cH^{p}}, \cD_{\cH} $.
They will be denoted by
$ K_{\cH^{p}}(X), K_{\cH}(X) $ if
$ Y $ is empty.

We define Deligne cohomology and absolute Hodge cohomology in the sense of
Beilinson ([3], [4]) by
$$
\aligned
{H}_{\cD}^{i}(X,Y;\bZ(j))
&= H^{i}\Gamma_{\cD}(K_{\cH}(X,Y)(j)),
\\
{H}_{\AH}^{i}(X,Y;\bZ(j))
&= H^{i}\Gamma_{\cH}(K_{\cH}(X,Y)(j)),
\endaligned
$$
See also [11], [12], [13], [19], etc.
We will omit
$ Y $ if it is empty.

By definition we have a natural morphism
$$
{H}_{\AH}^{i}(X,Y;\bZ(j)) \to
{H}_{\cD}^{i}(X,Y;\bZ(j)).
\leqno(4.2.3)
$$
\vglue12pt

 4.3. {\it Short exact sequences.} Since higher extensions vanish in
$ \MHS $, every complex is represented by a complex with zero differential
in
$ D^{b}\MHS $.
We see that
$ K \in \cD_{\cH} $ corresponds by (4.1.4) (noncanonically) to
$$
\mathbold{\oplus}_{i} (H^{i}K)[-i] \in D^{b}\MHS.
\leqno(4.3.1)
$$
Then, using the
$ t $-structure on
$ \cD_{\cH} $, we have a canonical exact sequence
$$
0 \to {\Ext}_{\MHS}^{1}(\bZ,H^{i-1}K(j)) \to H^{i}
\Gamma_{\cH}(K(j)) \to \Hom_{\MHS}(\bZ,H^{i}K(j))
\to 0.
\leqno(4.3.2)
$$

Similarly, we have
$$
0 \to J(H^{i-1}K(j)) \to H^{i}\Gamma_{\cD}(K (j))
\to H^{i}K_{\bZ}(j) \cap F^{j}H^{i}K_{\bC}
\to 0,
\leqno(4.3.3)
$$
where we put
$$
\aligned
J(H(j))
&= H_{\bC}/(H_{\bZ}(j) + F^{j}H_{\bC}),
\\
H_{\bZ}(j) \cap F^{j}H_{\bC}
&= \Ker(H_{\bZ}(j)
\to H_{\bC}/F^{j}H_{\bC}),
\endaligned
$$
for a mixed Hodge structure
$ H = (H_{\bZ}, (H_{\bQ},W), (H_{\bC};F,W)) $.

Comparing (4.3.2--3), we see that
$$
{H}_{\cD}^{i}(X_{\ssb},\bZ(j)) = {H}_{\AH}^{i}(X_{\ssb},
\bZ(j))
\leqno(4.3.4)
$$
if
$ H^{i-1}(X_{\ssb},\bZ) $ and
$ H^{i}(X_{\ssb},\bZ) $ have weights
$ \le 2j $ (using [8]).

\vglue12pt  4.4.\enspace {\it Weight spectral sequence.} Let
$ W $ be the weight filtration of\break
$ K_{\cH^{p}}(\osX\asD) $
in
$ \cC_{\cH^{p}} $, and
$ {\tD}_{j}^{k} $ the disjoint union of the intersections of
$ k $ irreducible components of
$ D_{j} $.
See [10].
By definition we have
$$
{\Gr}_{-p}^{W}K_{\cH^{p}}(\osX\asD)=
\mathbold{\oplus}_{k\ge 0} K_{\cH^{p}}({\tD}_{p+k}^{k})(-k)[-p-2k],
\leqno(4.4.1)
$$
and this gives the integral weight spectral sequence (2.3.2).
It depends on the choice of the compactification
$ \osX $ of
$ X_{\ssb} $.
By {\it loc.\ cit.}\ this spectral sequence degenerates at
$ E_{2} $ modulo torsion.

Let
$ {}^{o}W_{j}X_{\ssb} $ be the cofiltration in (3.1).
Then
$$
{\Gr}_{-p}^{W}K_{\cH^{p}}({}^{o}W^{j}\osX\asD) =
\mathbold{\oplus}_{k\ge 0} K_{\cH^{p}}
({}^{o}W^{j}{\tD}_{p+k}^{k})(-k)[-p-2k],
$$
where
$ {}^{o}W_{j}{\tD}_{p+k}^{k} = {\tD}_{p+k}^{k} $ if
$ p + k > j $ or
$ p - j = k = 0 $, and
$ {}^{o}W_{j}{\tD}_{p+k}^{k} = \emptyset $ otherwise.
This implies $\phantom{\sum^{|}}$
$$
H^{r}(W^{j}K_{\cH^{p}}(\osX\asD))=
H^{r}({X}_{\ssb}^{j}))\quad
\text{for }r \le j + 2. 
\leqno(4.4.2)
$$
Indeed, we have
$ W^{j}K_{\cH^{p}}({}^{o}W^{j}\osX\asD)=
W^{j}K_{\cH^{p}}(\osX\asD) $ and
$$
H^{r}(K_{\cH^{p}}({}^{o}W^{j}
\osX\asD)/
W^{j}K_{\cH^{p}}({}^{o}W^{j}
\osX\asD))
= 0\quad \text{for }r \le j + 2,
$$
where
$ W^{j} = W_{-j} $.
Note that the weight filtration on
$ H^{r}(X_{\ssb},\bZ) $ is shifted by
$ r $ as in [10] (i.e., it is induced by
$ \Dec\, W) $.

\vglue12pt  4.5. {\it Remark.} Assume
$ X $ is smooth proper.
Then
$ {H}_{\cD}^{i}(X,\bZ(j)) $ is the hypercohomology of the complex
$$
\bZ_{X^{\an}}(j) \to \cO_{X^{\an}} \to
{\Omega}_{{X}^{\an}}^{1} \to \dots \to
{\Omega}_{{X}^{\an}}^{j-1},
$$
where the degree of
$ \bZ_{X^{\an}}(j) $ is zero.
In particular, using the exponential sequence, we have for
$ j = 1 $
$$
{H}_{\cD}^{i}(X,\bZ(1))
= H^{i-1}(X^{\an}, \cO_{{X}^{\an}}^{*}).
\leqno(4.5.1)
$$
By (4.1.3) and (4.3.2), we get in the smooth case
$$
\aligned
&{H}_{\cD}^{i}(X,\bZ(1)) = {H}_{\AH}^{i}(X,\bZ(1))\quad
\text{for }i \le 2,
\\
&{H}_{\AH}^{i}(X,\bZ(1))\,\,\, \text{is torsion for }i > 3.
\endaligned
\leqno(4.5.2)
$$
Note that the algebraic cohomology
$ H^{i-1}(X, \cO_{X}^{*}) $ coincides with\break
$ H^{i-1}(X^{\an}, \cO_{{X}^{\an}}^{*}) $ for
$ i \le 2 $ by GAGA, and vanishes otherwise by (2.5.3).\break
In particular,
$ {H}_{\AH}^{i}(X,\bZ(1)) $ coincides with the algebraic cohomology\break
$ H^{i-1}(X, \cO_{X}^{*}) $ up to torsion for
$ i \ne 2 $ in the smooth case.
(If we consider a polarizable version, we get isomorphisms up to torsion
for every
$ i $.)

\section{Comparison}

In this section we prove Theorems (0.1--2).

\nonumproclaim{{5.1.~Theorem}} Let
$ X_{\ssb} $ and
$ \osX $ be as in {\rm (3.1)} with
$ k = \bC $.
Then there exist a decreasing filtration
$ W'' $ in a generalized sense on
$ \Gamma_{\cD}(K_{\cH}(X_{\ssb})(1)) $ and a canonical
morphism
$$
\Gamma (\osX^{\an},\cC^{\ssb}
({j}_{*}^{m}\cO_{{X}^{\an}}^{*} [-1])) \to
\Gamma_{\cD}(K_{\cH}(X_{\ssb})(1)),
\leqno(5.1.1)
$$
preserving the filtrations
$ W' $ and
$ W'' ${\rm ,} where
$ \cC^{\ssb} $ denotes the canonical flasque resolution{\rm ,}
$ {j}_{*}^{m}\cO_{{X}^{\an}}^{*} $ is the meromorphic extension{\rm ,} and
$ W' $ on the target is induced by the degree of
$ X_{\ssb} $ as in {\rm (3.1).}
Furthermore {\rm (5.1.1)} becomes a bifiltered quasi\/{\rm -}\/isomorphism by
replacing the target with
$ W^{\prime \prime -1} ${\rm ,} and the mapping cone of
$ {\Gr}_{W'}^{p} $ of {\rm (5.1.1)} is quasi\/{\rm -}\/isomorphic to
$ (\tau_{>1}j_{*}\cC^{\ssb}(\bZ_{{X}_{p}^{\an}}(1)))[-p] $. \pagebreak
\endproclaim

{\it Remark.} The filtration
$ W'' $ is not given by subcomplexes, but by morphisms of complexes
$ 0 = W^{\prime \prime 1} \to W^{\prime \prime 0} \to
W^{\prime \prime -1} \to
W^{\prime \prime -2} = \Gamma_{\cD}(K_{\cH}(X_{\ssb})(1)) $
compatible with
$ W' $.
See [5] and [27, 1.3].
For this we can define naturally the notion of bifiltered
quasi-isomorphism.

\demo{{P}roof of  {\rm (5.1)}}
Let
$ X $ be a smooth complex algebraic variety, and
$ \oX $ a smooth compactification such that
$ D := \oX \setminus X $ is a divisor with simple normal
crossings.
Then we have a short exact sequence
$$
0 \to \cO_{{\oX}^{\an}}^{*} \to
{j}_{*}^{m}\cO_{{X}^{\an}}^{*} \to
\mathbold{\oplus}_{i} \bZ_{{D}_{i}^{\an}} \to 0,
\leqno(5.1.2)
$$
where the
$ D_{i} $ are irreducible components of
$ D $.
This implies a canonical isomorphism in the derived category
$$
[\cO_{{\oX}^{\an}} \overset{\exp}\to\to
{j}_{*}^{m}\cO_{{X}^{\an}}^{*}]
= \tau_{\le 1}j_{*}\cC^{\ssb}(\bZ_{X^{\an}}(1)).
\leqno(5.1.3)
$$

We can verify that
$ \Gamma_{\cD}(K_{\cH}(X)(1)) $ is naturally quasi-isomorphic
to the complex of global sections of the shifted mapping cone
$$
[{j}_{*}\cC^{\ssb}(\bZ_{X^{\an}}(1))\oplus
\cC^{\ssb}(\sigma_{\ge 1}
{\Omega}_{{\oX}^{\an}}^{\ssb} (\log D)) \to
j_{*}\cC^{\ssb}({\Omega}_{{X}^{\an}}^{\ssb})],
$$
where
$ \cC^{\ssb} $ denotes the canonical flasque resolution of
Godement.
Since we have a natural morphism
$$
[\cO_{{X}^{\an}} \overset{\exp}\to\to \cO_{{X}^{\an}}^{*}]
\to {\Omega}_{{X}^{\an}}^{\ssb}
$$
using
$ d \log $ (see [10]), the above shifted mapping cone is
quasi-isomorphic to
$$
\bZ(1)_{\cD,X\aD} := [{j}_{*}
\cC^{\ssb}([\cO_{{X}^{\an}} \overset{\exp}\to\to
\cO_{{X}^{\an}}^{*}])\oplus \cC^{\ssb}(\sigma_{\ge 1}
{\Omega}_{{\oX}^{\an}}^{\ssb} (\log D)) \to
j_{*}\cC^{\ssb}({\Omega}_{{X}^{\an}}^{\ssb})].
$$

Consider then the shifted mapping cone
$$
\bZ(1)'_{\cD,X\aD} := [\cC^{\ssb}
([\cO_{{\oX}^{\an}}
\overset{\exp}\to\to {j}_{*}^{m}\cO_{{X}^{\an}}^{*}])\oplus
\cC^{\ssb}(\sigma_{\ge 1}
{\Omega}_{{\oX}^{\an}}^{\ssb}
(\log D)) \to j_{*}\cC^{\ssb}({\Omega}
_{{X}^{\an}}^{\ssb})].
$$
It has a natural morphism to
$ \bZ(1)_{\cD,X\aD} $, and defines
$ W^{\prime \prime -1} $.
Here we can replace
$ j_{*}\cC^{\ssb}({\Omega}_{{X}^{\an}}^{\ssb}) $ with
$ \cC^{\ssb}({\Omega}_{{\oX}^{\an}}^{\ssb}
(\log D)) $ because the image of
$ d \log $ is a logarithmic form.
So we see that
$ \bZ(1)'_{\cD,X\aD} $ is naturally
quasi-isomorphic to
$ \cC^{\ssb}({j}_{*}^{m}\cO_{{X}^{\an}}^{*})[-1] $.
Furthermore,
$ \bZ(1)'_{\cD,X\aD} $ has a subcomplex
$$
\bZ(1)_{\cD,\oX} := [\cC^{\ssb}
([\cO_{{\oX}^{\an}}
\overset{\exp}\to\to \cO_{{\oX}^{\an}}^{*}])
\oplus \cC^{\ssb}(\sigma_{\ge 1}
{\Omega}_{{\oX}^{\an}}^{\ssb}) \to
\cC^{\ssb}({\Omega}_{{\oX}^{\an}}^{\ssb})]
$$
which is naturally quasi-isomorphic to
$ \Gamma_{\cD}(K_{\cH}(\oX)(1)) $ and
$ \cC^{\ssb}(\cO_{{\oX}^{\an}}^{*})[-1] $.
This defines
$ W^{\prime \prime 0} $.

Thus we get a canonical filtered morphism
$$
\Gamma (\oX^{\an},\cC^{\ssb}
(j_{*}\cO_{{X}^{\an}}^{*}[-1])) \to
\Gamma_{\cD}(K_{\cH}(X)(1))
\leqno(5.1.4)
$$
whose mapping cone is isomorphic to
$ \tau_{>1}j_{*}\cC^{\ssb}(\bZ_{X^{\an}}(1)) $.

We apply this constriction to each component
$ \oX_{p} $ of
$ \osX $.
Then we get the morphism (5.1.1) preserving the filtrations
$ W' $ and
$ W'' $.
\enddemo

\nonumproclaim{{5.2.~Corollary}} For
$ X_{\ssb} $ as in {\rm (5.1),} we have a canonical morphism
$$
\Pic(X_{\ssb}) \to {H}_{\cD}^{2}(X_{\ssb},\bZ(1)),
\leqno(5.2.1)
$$
induced by {\rm (5.1.1).}
{\rm (}See {\rm (2.4)} for
$ \Pic(X_{\ssb}).) $
This is injective if
$ X_{j} $ is empty for
$ j < 0 ${\rm ,} and is bijective if furthermore
$ X_{0} $ is proper.
In particular{\rm ,} we have in the notation of {\rm (3.1)} and {\rm (4.1)}
$$
P_{\ge r}(\osX\asD) =
{H}_{\cD}^{r+2}({}^{o}W^{r}X_{\ssb},\bZ(1)).
\leqno(5.2.2)
$$
\endproclaim

\demo{{P}roof} This follows from (5.1) together with
(2.5.3).
\enddemo

{\it  Remark}. By (4.3.2) and (4.3.4) we get
$$
{H}_{\AH}^{r+i}({}^{o}W^{r}\osX,\bZ(1)) =
{H}_{\cD}^{r+i}({}^{o}W^{r}X_{\ssb},\bZ(1))
\quad \text{for }i \le 2,
\leqno(5.2.3)
$$
and a short exact sequence
$$
\aligned
0
&\to {\Ext}_{\MHS}^{1}(\bZ,H^{r+1}({}^{o}W^{r}X_{\ssb},
\bZ)(1)) \to
{H}_{\cD}^{r+2}({}^{o}W^{r}X_{\ssb},\bZ(1))
\\
&\to \Hom_{\MHS}(\bZ,H^{r+2}({}^{o}W^{r}X_{\ssb},
\bZ)(1)) \to 0,
\endaligned
\leqno(5.2.4)
$$
because
$ H^{r+i}({}^{o}W^{r}X_{\ssb},\bZ) $ has weights
$ \le i $ for
$ i \le 2 $ by (4.4.2).
Furthermore, using the spectral sequence (2.3.2), we see
$$
{\Ext}_{\MHS}^{1}(\bZ,H^{r+1}({}^{o}W^{r}\osX,
\bZ)(1)) = {\Ext}_{\MHS}^{1}(\bZ,H^{r+1}({}^{o}W^{r}X_{\ssb},
\bZ)(1)),
\leqno(5.2.5)
$$
and they have naturally a structure of semiabelian variety as is
well-known.
See [8], [10].
This coincides with the structure of semiabelian variety on
$ \cP_{\ge r}(\osX)^{0} $ by the following:

\nonumproclaim{{5.3.~Theorem}} Let
$ S $ be a smooth complex algebraic variety{\rm ,} and
$ (L,\gamma) \in P_{\ge r}(\osX\times S) $ as in
Remark after {\rm (3.1),} where
$ D_{\ssb} $ is empty.
Then the set\/{\rm -}\/theoretic map
$$
S(\bC) \to {H}_{\cD}^{r+2}({}^{o}W^{r}\osX,\bZ(1))
$$
defined by restricting
$ (L,\gamma) $ to the fiber at
$ s \in S(\bC) $ comes from a morphism of varieties.
\endproclaim

\demo{{P}roof} Let
$ \xi \in {H}_{\cD}^{r+2}({}^{o}W^{r}\osX\times
S,\bZ (1)) $ corresponding to
$ (L,\gamma) $ by the injective morphism
$$
P_{\ge r}(\osX\times S) \to {H}_{\cD}^{r+2}
({}^{o}W^{r}\osX\times S,\bZ(1))
\leqno(5.3.1)
$$
given by (5.2.1) for
$ \osX\times S $.
Since this morphism is compatible with the restriction to
$ X_{\ssb}\times \{s\} $, it is enough to show that the map
$$
s \mapsto \xi_{s} \in {H}_{\cD}^{r+2}({}^{o}W^{r}
\osX, \bZ(1)) = {H}_{\AH}^{r+2}({}^{o}W^{r}
\osX,\bZ(1))
$$
is algebraic, where the last isomorphism follows from (5.2.3).

We first replace the Deligne cohomology of
$ {}^{o}W^{r}\osX\times S $ with the absolute Hodge
cohomology.
Since the restriction to the fiber is defined at the level of mixed Hodge
complexes, it is compatible with the canonical morphism induced by
$ \Gamma_{\cH} \to \Gamma_{\cD} $.
So it is enough to show that
$ \xi $ belongs to the image of the absolute Hodge cohomology.
But this is verified by using
$ {}^{o}W^{r}(\osX\times S) $ which is defined by
replacing the
$ r $-th component of
$ {}^{o}W^{r}\osX\times S $ with
$ {}^{o}W^{r}(\oX_{r}\times \oS) $, where
$ \oS $ is a smooth compactification of
$ S $.
Indeed, the Deligne cohomology and the absolute Hodge cohomology
coincide for this by (4.3.4), and the line bundle
$ L $ can be extended to
$ \oX_{r}\times \oS $.

Now we reduce the assertion to the case
$$
\xi_{s} \in {\Ext}_{\MHS}^{1}(\bZ,H^{r+1}({}^{o}W^{r}
\osX, \bZ)(1)).
$$
Indeed, the image of
$ \xi_{s} $ in
$ \Hom_{\MHS}(\bZ,H^{r+2}({}^{o}W^{r}X_{\ssb},\bZ)(1)) $ by
(5.2.4) is constant, and we may assume it is zero by adding the pull-back
of an element of
$ {H}_{\cD}^{r+2}({}^{o}W^{r}\osX,\bZ(1)) $ by
the projection
$ \osX\times S \to \osX $.

We then claim that
$ \{\xi_{s}\}_{s\in S(\bC)} $ is an admissible normal function in
the sense of [29], i.e., it defines an extension between constant
variations of mixed Hodge structures on
$ S $ and the obtained extension is an admissible variation of mixed Hodge
structure in the sense of Steenbrink-Zucker [32] and Kashiwara [20].
(Actually it is enough to show that
$ \{\xi_{s}\} $ is an analytic section for the proof of (5.4),
because an analytic structure of a semiabelian variety is equivalent
to an algebraic structure [10].)

Choosing a splitting of the exact sequence (4.3.2) for
$ K_{\cH}({}^{o}W^{r}\osX\times S) $, we get a
decomposition of
$ \xi $:
$$
\aligned
\xi'
&\in {\Ext}_{\MHS}^{1}(\bZ,H^{r+1}({}^{o}W^{r}\osX
\times S,\bZ)(1)),
\\
\xi''
&\in \Hom_{\MHS}(\bZ,H^{r+2}({}^{o}W^{r}\osX
\times S,\bZ)(1)).
\endaligned
$$
Then only the following K\"unneth components contribute to the restriction
to the fiber at
$ s $:
$$
\aligned
\xi^{\prime 0}
&\in {\Ext}_{\MHS}^{1}(\bZ,
H^{r+1}({}^{o}W^{r}\osX,\bZ)(1)
\otimes H^{0}(S,\bZ)),
\\
\xi^{\prime \prime 1}
&\in \Hom_{\MHS}(\bZ,
H^{r+1}({}^{o}W^{r}\osX,\bZ)(1)
\otimes H^{1}(S,\bZ)).
\endaligned
$$
Clearly,
$ \xi^{\prime 0} $ gives a constant section (where we may assume
$ S $ connected), and the restriction of
$ \xi^{\prime \prime 1} $ is well-defined modulo constant section.
We will show that the restrictions of
$ \xi^{\prime \prime 1} $ to the points of
$ S $ form an admissible normal function.

The restriction is also defined by applying the functor
$ \Gamma_{\cH} $ to the restriction morphism of mixed Hodge
complexes
$$
K_{\cH}({}^{o}W^{r}\osX\times S)(1) \to
K_{\cH}({}^{o}W^{r}\osX)(1).
\leqno(5.3.2)
$$
By (4.1.4) this corresponds to the tensor of
$ \mathbold{\oplus}_{i} H^{i}({}^{o}W^{r}\osX,\bZ)[-i ] $
with the restriction morphism under the inclusion
$ \{s\} \to S $:
$$
{\bf R}\Gamma(S,\bZ)\,(\simeq
\mathbold{\oplus}_{i} H^{i}(S,\bZ)[-i]) \to \bZ\quad
\text{in }D^{b}\MHS.
\leqno(5.3.3)
$$
Note that the restriction of the morphism to
$ H^{i}(S,\bZ)[-i] $ vanishes for
$ i > 1 $, and choosing
$ s_{0} \in S $, the restriction to
$ H^{1}(S,\bZ)[-1] $ for
$ s \ne s_{0} $ is expressed by the short exact sequence
$$
0 \to \bZ \to H^{1}(S,\{s_{0},s\};\bZ)
\to H^{1}(S,\{s_{0}\};\bZ) \,(= H^{1}(S,\bZ)) \to 0,
$$
using the corresponding distinguished triangle.
Here
$ H^{1}(S,\bZ)[-1] \to \bZ $ is zero for
$ s = s_{0} $.

The restriction of
$ \xi^{\prime 0} $ and
$ \xi^{\prime \prime 1} $ is then obtained by tensoring (5.3.3) with
$ H := H^{r+1}({}^{o}W^{r}\osX,\bZ)(1) $, and
applying the functor
$$
{\bf R}\Hom_{D^{b}\MHS}(\bZ[-1],*).
$$
For
$ \xi^{\prime \prime 1} $, it is given by taking the pull-back of
the above short exact sequence tensored with
$ H $ by
$ \xi^{\prime \prime 1} $.
Then we can construct the extended variation of mixed Hodge structure by
using the diagonal of
$ S\times S $ as in [28, 3.8].
This shows that
$ \{\xi_{s}\} $ determines an admissible normal function which will be
denoted by
$ \rho $.

Let
$ G $ be the semiabelian variety defined by
$${\Ext}_{\MHS}^{1}(\bZ,
H^{r+1}({}^{o}W^{r}\osX,\bZ)(1))$$
(see [8], [10]).
Then
$ \rho $ is a holomorphic section of
$ G^{\an}\times S^{\an} \to S^{\an} $.
We have to show that this is algebraic using the property of admissible
normal functions.
Since
$ G $ is semiabelian, there exist a torus
$ T $ and an abelian variety
$ A $ together with a short exact sequence
$$
0 \to T \to G \to A \to 0.
$$
As a variety,
$ G $ may be viewed as a principal
$ T $-bundle.
We choose an isomorphism
$$
T = (\bG_{m})^{n}.
$$
This gives compactifications
$ \oT = (\bP^{1})^{n} $ of
$ T $, and also
$ \oG $ of
$ G $.

Then, by GAGA, it is enough to show that the admissible normal function
$ \rho $ is extended to a holomorphic section of
$ \oG^{\an}\times \oS^{\an} \to
\oS^{\an} $, where
$ \oS $ is an appropriate smooth compactification of
$ S $ such that
$ \oS \setminus S $ is a divisor with normal crossings.
Here we may assume
$ n = 1 $ using the projections
$ \oT \to \bP^{1} $, because
$ \oG $ is the fiber product of the
$ \bP^{1} $-bundles over
$ A $.

So the assertion follows from the same argument as in [29, 4.4].
Indeed, the group of connected components of the fiber of the N\'eron
model of
$ G\times S $ at a generic point of
$ \oS \setminus S $ is isomorphic to
$ \bZ $ by an argument similar to (2.5.5) in {\it loc.\ cit.}, and this
corresponds to the order of zero or pole in an appropriate sense
of a local section.
By blowing up further, we may assume that these orders along any
intersecting two of the irreducible components of the divisor have the
same sign (including the case where one of them is zero).
Then it \pagebreak can be extended to a section of
$ \oG^{\an}\times \oS^{\an} \to
\oS^{\an} $ as in the case of meromorphic functions
(which corresponds to the case
$ A = 0) $.
This finishes the proof of (5.3).
\enddemo

{\it Remark}. It is easy to show that the map
$ S(\bC) \to {H}_{\cD}^{r+2}({}^{o}W^{r}\osX,\bZ(1)) $ is analytic
using the long exact sequence associated with the direct image of the
distinguished triangle
$ \to \bZ(1) \to \cO \to \cO^{*} \to $ under the projection
$ \oX_{\ssb}^{\an}\times S^{\an} \to S^{\an} $.
This implies that the natural algebraic structure on
$ {H}_{\cD}^{r+2}({}^{o}W^{r}\osX,\bZ(1)) $ is compatible with the
one obtained by (3.2) and (5.2.2).

\nonumproclaim{{5.4.~Theorem}} With the notation of {\rm (4.2),} let
$ W $ denote the decreasing weight filtration on the mixed Hodge complex
$ K := K_{\cH^{p}}(\osX\asD) $ {\rm (}\/i.e.\ 
$ W^{j} = W_{-j}) $.
For a mixed Hodge structure
$ H ${\rm ,} let
$ H_{(1)} $ denote the maximal mixed Hodge structure contained in
$ H $ and such that
$ {\Gr}_{F}^{p} = 0 $ for
$ p \notin \{0,1\} $.
Then we have natural isomorphisms of mixed Hodge structures
$$
\align
&\tag 5.4.1
\\ \noalign{\vskip-18pt}
r_{\cH}(G_{r}(\osX\asD))(-1)
&= W_{1}H^{r}(W^{r-2}K)_{\fr},
\\ \navs
&\tag 5.4.2  \\ \noalign{\vskip-18pt}
r_{\cH}(\Gamma'_{r}(\osX\asD))(-1)
&= H^{r}(W^{r-2}K)_{(1)}/W_{1}H^{r}(W^{r-2}K),
\\ \navs
\tag 5.4.3 \\ \noalign{\vskip-18pt}
r_{\cH}(M'_{r}(\osX\asD))(-1)
&= H^{r}(W^{r-2}K)_{(1)}/W_{1}H^{r}(W^{r-2}K)_{\tor},
\endalign
$$
and surjective morphisms of mixed Hodge structures
$$
\gather
r_{\cH}(\Gamma_{r}(\osX\asD))(-1)
\to W_{2}H^{r}(W^{r-3}K)_{(1)}/W_{1}H^{r}(W^{r-2}K),
\tag 5.4.4
\\ \navs
r_{\cH}(M_{r}(\osX\asD))(-1)
\to W_{2}H^{r}(W^{r-3}K)_{(1)}/W_{1}H^{r}(W^{r-3}K)_{\tor},
\tag 5.4.5
\endgather
$$
whose kernels are torsion{\rm ,} and vanish if
$ H^{2}(\oX_{r-3},\bZ) $ is of type
$ (1,1) $.
\endproclaim

\demo{{P}roof} By (4.4.2) and (5.3) we have
$$
r_{\cH}(\cP_{\ge r-1}(\osX)^{0})(-1) =
H^{r}(W^{r-1}K)_{\fr}.
$$
Since
$ r_{\cH}(\cP_{r-2}(\osX)^{0})(-1) =
H^{1}(\oX_{r-2},\bZ)_{\fr} = H^{r-1}
({\Gr}_{W}^{r-2}K)_{\fr} $, we get
$$
r_{\cH}(G_{r}(\osX\asD))(-1) =
\Coker(\partial : H^{1}(\oX_{r-2},\bZ) \to
H^{r}(W^{r-1}K))_{\fr},
$$
using the right exactness of
$ \Ext^{1}(\bZ,*) $.
So (5.4.1) follows from the long exact sequence of mixed Hodge structures
$$
\aligned
H^{1}(\oX_{r-2},\bZ)
&\to H^{r}(W^{r-1}K)\to H^{r}(W^{r-2}K)
\\
&\to H^{2}(\oX_{r-2},\bZ)
 \oplus \Gamma (\tD_{r-1},\bZ)(-1)
\overset{\partial}\to\to H^{r+1}(W^{r-1}K).
\endaligned
\leqno(5.4.6)
$$
Here
$ H^{i+r}{\Gr}_{W}^{r-2}K = H^{i+2}(\oX_{r-2},
\bZ) \oplus H^{i}(\tD_{r-1},\bZ)(-1) $ for
$ i \le 1 $ by (4.4).
Note that
$ W_{i}H^{r}(W^{j}K) = \Im(H^{r}(W^{r-i}K) \to
H^{r}(W^{j}K) $ for
$ r - i \ge j $.

Similarly, we have by (5.2.2)
$$
P_{\ge r-1}(\osX\asD)/P_{\ge
r-1}(\osX)^{0} = \Hom_{\MHS}(\bZ,H^{r+1}
(W^{r-1}K)(1)).
$$
Furthermore, the right vertical morphism of (3.1.4) is identified with
the image by the functor
$ \Hom_{\MHS}(\bZ(-1),*) $ of the last morphism
$ \partial $ in (5.4.6).
This is verified by using a canonical morphism of triangles
in
$ \cD_{\cH} $
$$
\matrix
 \lrar&  W^{r-1}K &\lrar&  W^{r-2}K &\lrar&  {\Gr}_{W}^{r-2}K &\lrar 
\\
\navs
 &\Big\downarrow&&\Big\downarrow&&\Big\downarrow
\\ \navs
 \lrar&  K_{\cH}({}^{o}W^{r-1}X_{\ssb}) &\lrar&  K_{\cH}
({}^{o}W^{r-2}X_{\ssb}) &\lrar&  K_{\cH}
({\Gr}_{{}^{o}W}^{r-2}X_{\ssb}) &\lrar 
\endmatrix
$$
where
$ {\Gr}_{{}^{o}W}^{r-2}X_{\ssb} $ is defined by the shifted mapping
cone.
Note that the left part of the diagram commutes without homotopy so that
the morphism of the mapping cones is canonically defined, and the
filtration
$ \Dec\, W $ on
$ W^{i}K/W^{j}K $ coincides with the filtration induced by
$ \Dec\, W $ on
$ K $.

So (5.4.2) follows from (5.4.6), because the functor
$ \Hom_{\MHS}(\bZ(-1),*) $ is left exact, and its restriction to pure
Hodge structures of weight
$ 2 $ is identified with the functor
$ H \to H_{(1)} $.
Note that a pure Hodge structure of type
$ (1,1) $ is identified with a
$ \bZ $-module.

For (5.4.3), we consider the extension class in
$$
{\Ext}_{\MHS}^{1}(\Gamma'_{r}(\osX\asD),
W_{1}H^{r}(W^{r-2}K)_{\fr}(1)),
$$
which is induced by (5.4.2) together with the natural exact sequence.
Here
$ \Gamma'_{r}(\osX\asD) $ is
identified with a mixed Hodge structure of type (0,0).
In particular, the extension class is equivalent to the induced morphism
$$
\Hom_{\MHS}(\bZ,\Gamma'_{r}(\osX\asD)) \to
{\Ext}_{\MHS}^{1}(\bZ,W_{1}H^{r}(W^{r-2}K)_{\fr}(1)).
\leqno(5.4.7)
$$
So we have to show for the proof of (5.4.3) that this morphism is
identified by (5.4.1) with
$$
\Gamma'_{r}(\osX\asD)
\to G_{r}(\osX\asD).
$$

By (5.2) and (4.4.2), the commutative diagram (3.1.4) is identified with a
morphism of the short exact sequences (4.3.2) induced by the morphism of
absolute Hodge cohomologies
$$
\partial : H^{r}\Gamma_{\cH}({\Gr}_{W}^{r-2}K(1)) \to
H^{r+1}\Gamma_{\cH}(W^{r-1}K(1)).
$$
The last morphism
$ \partial $ is induced by the ``boundary map'' of the first
distinguished triangle of the above diagram (which is defined by using the
mapping cone).
Then (5.4.3) follows from (4.1.3--4) by using (5.5) below,
because
$ H^{i}\Gamma_{\cH} $ is identified with
$ {\Ext}_{\MHS}^{i}(\bZ,*) $ by the equivalence of categories (4.1.4)
due to (4.1.3).

Finally, the surjectivity of (5.4.4) and (5.4.5) follows from the exact
sequence
$$
H^{2}(\oX_{r-3},\bZ) \oplus H^{0}(\tD_{r-2},
\bZ)(-1)\to H^{r}(W^{r-2}K)
\overset{\partial}\to\to H^{r}(W^{r-3}K),
\leqno(5.4.8)
$$
by comparing the (1,1) part of
$ \Coker\,{\Gr}_{2}^{W}\partial $ with the cokernel of the (1,1) part of
$ {\Gr}_{2}^{W}\partial $ (where
$ {\Gr}_{2}^{W}\Coker\,\partial = \Coker\,{\Gr}_{2}^{W}\partial $
because
$ H^{r}(W^{r-2}K) = W_{2}H^{r}(W^{r-2}K) $).
The kernels of (5.4.4) and (5.4.5) come from the difference between these,
and vanish if
$ H^{2}(\oX_{r-3},\bZ) $ is of type (1,1).
Note that the intersection of
$ W_{1}H^{r}(W^{r-2}K) $ with
$$
\Im(H^{2}(\oX_{r-3},\bZ) \oplus H^{0}(\tD_{r-2},
\bZ)(-1) \to H^{r}(W^{r-2}K))
$$
is torsion by the strict compatibility of the weight filtration.
This completes the proof of (5.4).
\enddemo

 5.5. {\it Remark.} Let
$ \cA $ be an abelian category such that
$ \Ext^{i}(A,B) = 0 $ for any objects
$ A, B $ and
$ i > 1 $.
Let
$ A_{0} \in \cA $, and
$ F(K) = {\bf R}\Hom(A_{0},K) $ for
$ K \in D^{b}\cA $ so that we have a canonical short exact sequence
$$
0 \to \Ext^{1}(A_{0},H^{-1}K) \to F(K) \to
\Hom(A_{0},H^{0}K) \to 0.
$$
Let
$ \to K' \to K \to K'' \to $ be a
distinguished triangle in
$ D^{b}\cA $ with
$ \partial : K'' \to K'[1] $ the boundary map.
Consider the morphism
$$
\Hom(A_{0},\Ker\,H^{0}\partial) \to \Ext^{1}(A_{0},
\Coker\,H^{-1}\partial),
\leqno(5.5.1)
$$
obtained by the snake lemma together with the right exactness of
$ \Ext^{1}(A_{0},*) $.
(Here
$ H^{i}\partial $ is the abbreviation of
$ H^{i}\partial : H^{i}K'' \to H^{i+1}K' $.)
Then (5.5.1)
coincides with the morphism induced by the short exact sequence
$$
0 \to \Coker\,H^{-1}\partial \to H^{0}K \to
\Ker\,H^{0}\partial \to 0.
\leqno(5.5.2)
$$

Indeed,
$ K', K'' $ are represented by complexes with {\it zero} differential, and
the assertion is reduced to the case where
$ K'' = A, K' = B $ with
$ A, B \in \cA $ (considering certain subquotients of
$ K', K'') $.
Then it follows from the well-known bijection between the extension group
in the derived category and the set of extension classes in the usual
sense.

\vglue12pt {5.6. {\it Proof of} (0.1--3).}
We take a resolution of
$ [Y \to X] $ so that the associated integral weight filtration is
defined independently of the choice of the resolution as in the proof of
(2.3) (e.g. we can take the simplicial resolution of Gillet and Soul\'e
[14, 3.1.2] if
$ X $ is proper).
Then by (5.4), it is enough to show that the kernel of the canonical
morphism
$$
W_{2}H^{r}(W^{r-3}K) \to H^{r}(K)
$$
is torsion, and it is contained in
$ W_{1}H^{r}(W^{r-3}K)_{\tor} $ if
$ {E}_{2}^{p,r-1-p} = 0 $ for
$ p \le r - 4 $.
But these can be verified by using a natural morphism between the
weight spectral sequences (2.3.2) converging to
$ H^{\ssb}(W^{r-3}K) $ and
$ H^{\ssb}(K), $ \pagebreak because the spectral sequences degenerate at
$ E_{2} $ modulo torsion.
Note that the weight filtration on the cohomology is shifted by the
degree, and
$ W_{2} $ is induced by
$ W^{r-2} $.
This completes the proof of (0.1--3).

\references

[1]
\name{L.\ Barbieri-Viale},  On algebraic $1$-motives related to Hodge cycles,
in {\it Algebraic Geometry\/} - {\it A Volume in Memory of
Paolo Francia\/}, W.\  de Gruyter, New York, 2002, 25--60.

[2]
\name{L.\ Barbieri-Viale} and \name{V.\ Srinivas},  {\it Albanese and Picard
$1$-motives\/},
{\it M{\rm \'{\it e}}m.\ Soc.\ Math.\ France\/} {\bf 87}, Paris, 2001.

[3]
\name{A.\ Beilinson}, Higher regulators and values of
$L$-functions, {\it J.\ Soviet Math\/}.\ {\bf 30}  (1985), 2036--2070.

[4]
\bibline,  Notes on absolute Hodge cohomology, {\it Applications of
Algebraic $K$-theory to Algebraic Geometry and Number Theory\/}, {\it
Part\/} II
(Boulder, CO, 1983), {\it Contemp.\
Math\/}.\  {\bf 55} (1986), 35--68.

[5]
\name{A.\ Beilinson, J.\ Bernstein}, and \name{P.\ Deligne},  {\it Faisceaux
pervers\/},
 {\it Analysis and Topology on Singular Spaces\/}, I (Luminy, 1981), 
{\it Ast{\rm \'{\it e}}risque\/} {\bf 100}, 5--171, Soc.\ Math\. France, Paris, 1982.

[6]
\name{J.\ Biswas} and \name{V.\ Srinivas},  A Lefschetz $(1,1)$ theorem for normal
projective complex varieties, {\it Duke Math.\ J\/}.\ {\bf 101} (2000), 427--458.

[7]
\name{S.\ Bloch},  Algebraic cycles and higher
$K$-theory, {\it Adv.\ in Math\/}.\ {\bf 61}  (1986), 267--304.

[8]
\name{J.\ Carlson},  Extensions of mixed Hodge structures, in {\it Journ{\rm \'{\it e}}es de
G{\rm \'{\it e}}om{\rm \'{\it e}}trie Alg{\rm \'{\it e}}brique d\/{\rm '}\/Angers\/},  1979, 107--127, Sijthoff and
Noordhoff, Germantown, MD, 1980.

[9]
\bibline,  The one-motif of an algebraic surface, {\it Compositio
Math\/}.\ {\bf 56} 
(1985), 271--314.

[10]
\name{P.\ Deligne},  Th\'eorie de Hodge I, {\it Actes du Congr\`es
International des Math{\rm \'{\it e}}maticiens\/} (Nice, 1970) 
{\bf 1}, 425--430; II, {\it Publ.\ Math.\ IHES\/} {\bf 40}  (1971), 5--57; 
III, {\it ibid\/}.\ {\bf 44} 
(1974), 5--77.

[11]
\name{F.\ El Zein} and \name{S.\ Zucker},  Extendability of normal functions associated
to algebraic cycles, in {\it Topics in Transcendental Algebraic
Geometry\/} (Princeton, NJ, 1981/1982),  269--288, 
{\it Ann.\ of Math.\ Studies\/} {\bf 106},
 Princeton Univ.\ Press, Princeton, NJ, 1984.

[12]
\name{H.\ Esnault} and \name{E.\ Viehweg},  Deligne-Beilinson cohomology, in 
{\it Perspect.\ Math\/}.\ {\bf 4},  43--91, Academic Press, Boston, MA, 1988.

[13]
\name{H.\ Gillet}, Deligne homology and Abel-Jacobi maps, {\it Bull.\ Amer.\
Math.\ Soc\/}.\ 
{\bf 10}  (1984), 285--288.

[14]
\name{H.\ Gillet} and \name{C.\ Soul\'e},  Descent, motives and $K$-theory, {\it J.\
Reine Angew.\ 
Math\/}.\ {\bf 478}  (1996), 127--176.

[15]
\name{A.\ Grothendieck},  {\it Fondements de la G{\rm \'{\it e}}om{\rm \'{\it e}}trie Alg{\rm \'{\it e}}brique\/}, 
Benjamin, New York, 1966.

[16]
\name{F.\ Guill\'en} and \name{V.\ Navarro Aznar},  Un crit\`ere d'extension d'un
foncteur d\'efini sur les sch\'emas lisses, {\it Publ.\ Math.\ IHES} {\bf 95} (2002), 1--91.

[17]
\name{F.\ Guill\'en,  V.\ Navarro Aznar,  P.\ Pascual-Gainza}, and \name{F.\ Puerta}, 
{\it Hyperr{\rm \'{\it e}}solutions Cubiques et Descente Cohomologique\/},
{\it Lecture Notes in Math\/}.\ {\bf 1335}, Springer-Verlag, New York,
1988.

[18]
\name{H.\ Hironaka},  Flattening theorem in complex-analytic geometry,
{\it Amer.\ J.\ Math\/}.\ {\bf 97} (1975),  503--547.

[19]
\name{U.\ Jannsen},  Deligne homology, Hodge-$D$-conjecture, and motives,
in {\it Perspect.\ Math\/}.\ {\bf 4}, 305--372,
 Academic Press, Boston, MA, 1988.

[20]
\name{M.\ Kashiwara},  A study of variation of mixed Hodge structure, {\it
Publ.\ Res.\ Inst.\ Math.\ Sci\/}.\  {\bf 22} (1986), 991--1024.

[21]
\name{B.\ Mazur} and \name{W.\ Messing}, Universal extensions and one dimensional
crystalline cohomology, {\it Lecture Notes in Math\/}.\ {\bf 370},
Springer-Verlag, New York, 
 1974.

[22]
\name{D.\ Mumford},  {\it Abelian Varieties\/}, 
 Oxford Univ.\ Press, London, 1970.

[23]
\name{J.\ Murre},  On contravariant functors from the category of pre-schemes over
a field into the category of abelian groups (with an application to the
Picard functor), {\it Publ.\ Math.\ IHES\/} 
{\bf 23}   (1964),
5--43.

[24]
\name{F.\ Oort}, {\it Commutative Group Schemes\/}, {\it Lecture Notes in
Math\/}.\ {\bf 15}, 
Springer-Verlag, New York, 1966.

[25]
\name{N.\ Ramachandran}, Duality of Albanese and Picard $1$-motives, {\it
$K$-Theory\/} {\bf 22} (2001), 271--301.

[26]
\bibline,  One-motives and a conjecture of Deligne, preprint
http://arxiv.org/abs/ \break math.AG/9806117 (revised
version v3, February 2003); {\it J. Algebraic Geometry\/}, to appear.

[27]
\name{M.\ Saito},  Modules de Hodge polarisables, {\it Publ.\ RIMS 
Kyoto Univ\/}.\ {\bf 24}
(1988), 849--995.

[28]
\bibline,  Extension of mixed Hodge modules, {\it Compositio Math\/}.\
{\bf 74} (1990), 
209--234.

[29]
\bibline,  Admissible normal functions, {\it J.\ Algebraic Geom\/}.\
{\bf 5}  (1996),
235--276.

[30]
\bibline,  Mixed Hodge complexes on algebraic varieties, {\it Math.\
Ann\/}.\ {\bf 316}
(2000), 283 --331.

[31]
\bibline,  Bloch's conjecture, Deligne cohomology and higher Chow
groups, preprint\break (math.AG/9910113).

[32]
\name{J.\ Steenbrink} and \name{S.\ Zucker},  Variation of mixed Hodge structure.\ I, 
{\it Invent. Math\/}.\ {\bf 80} 
 (1985), 489--542.

[33]
\name{J.-L.\ Verdier},  Cat\'egories d\'eriv\'ees, in SGA 4 1/2, {\it Lecture
Notes in Math\/}.\
{\bf 569}, 262--311, Springer-Verlag, New York,  1977.
\endreferences

\bye

MEMOMAT

\noindent
Universit\`a di Roma ``La Sapienza''

\noindent
Via A. Scarpa, 16

\noindent
I--00161 Roma

\noindent
E-Mail: barbieri\@dmmm.uniroma1.it

\bigskip
\noindent
Andreas Rosenschon

\noindent
Department of Mathematics

\noindent
Duke University

\noindent
Durham, NC 27708

\noindent
U.S.A

\noindent
E-Mail: axr\@math.duke.edu

\bigskip
\noindent
Morihiko Saito

\noindent
RIMS Kyoto University,

\noindent
Kyoto 606--8502

\noindent
Japan

\noindent
E-Mail: msaito\@kurims.kyoto-u.ac.jp

\bigskip
\noindent

\bye